\documentclass[12pt]{article}
\usepackage[utf8x]{inputenc}
\usepackage{amssymb, amsmath, amsthm, amscd}
\usepackage[pdftex]{graphics}
\usepackage[all]{xy}
\usepackage{bbm}
\usepackage{enumitem}
\usepackage{setspace}
\usepackage[colorlinks=true]{hyperref}
\usepackage{centernot}
\hypersetup{colorlinks   = true}
\hypersetup{linkcolor=blue}
\usepackage{tikz}
\usetikzlibrary{cd}
\usetikzlibrary{matrix}

\begin{document}
\newtheorem{The}{Theorem}[section]
\newtheorem{Lem}[The]{Lemma}
\newtheorem{Prop}[The]{Proposition}
\newtheorem{Cor}[The]{Corollary}
\newtheorem{Rem}[The]{Remark}
\newtheorem{Obs}[The]{Observation}
\newtheorem{SConj}[The]{Standard Conjecture}
\newtheorem{Titre}[The]{\!\!\!\! }
\newtheorem{Conj}[The]{Conjecture}
\newtheorem{Question}[The]{Question}
\newtheorem{Prob}[The]{Problem}
\newtheorem{Def}[The]{Definition}
\newtheorem{Not}[The]{Notation}
\newtheorem{Claim}[The]{Claim}
\newtheorem{Conc}[The]{Conclusion}
\newtheorem{Ex}[The]{Example}
\newtheorem{Fact}[The]{Fact}
\newtheorem{Formula}[The]{Formula}
\newtheorem{Formulae}[The]{Formulae}
\newtheorem{The-Def}[The]{Theorem and Definition}
\newtheorem{Prop-Def}[The]{Proposition and Definition}
\newtheorem{Cor-Def}[The]{Corollary and Definition}
\newtheorem{Conc-Def}[The]{Conclusion and Definition}
\newtheorem{Lem-Def}[The]{Lemma and Definition}
\newtheorem{Terminology}[The]{Note on terminology}
\newtheorem{Construction}[The]{Construction}
\newcommand{\C}{\mathbb{C}}
\newcommand{\R}{\mathbb{R}}
\newcommand{\N}{\mathbb{N}}
\newcommand{\Z}{\mathbb{Z}}
\newcommand{\Q}{\mathbb{Q}}
\newcommand{\Proj}{\mathbb{P}}
\newcommand{\Rc}{\mathcal{R}}
\newcommand{\Oc}{\mathcal{O}}

\newcommand\slawek[1]{{\textcolor{red}{#1}}}
\newcommand\dan[1]{{\textcolor{blue}{#1}}}

\begin{center}

{\Large\bf A Generalised Volume Invariant for Aeppli Cohomology Classes of Hermitian-Symplectic Metrics}

\end{center}

\begin{center}
{\large S\l{}awomir Dinew and Dan Popovici}

\end{center}

\vspace{1ex}

\noindent{\small{\bf Abstract.} We investigate the class of compact complex Hermitian-symplectic manifolds $X$. For each Hermitian-symplectic metric $\omega$ on $X$, we introduce a functional acting on the metrics in the Aeppli cohomology class of $\omega$ and prove that its critical points (if any) must be K\"ahler when $X$ is $3$-dimensional. We go on to exhibit these critical points as maximisers of the volume of the metric in its Aeppli class and propose a Monge-Amp\`ere-type equation to study their existence. Our functional is further utilised to define a numerical invariant for any Aeppli cohomology class of Hermitian-symplectic metrics that generalises the volume of a K\"ahler class. We obtain two cohomological interpretations of this invariant.  Meanwhile, we construct an invariant in the form of an $E_2$-cohomology class, that we call the $E_2$-torsion class, associated with every Aeppli class of Hermitian-symplectic metrics and show that its vanishing is a necessary condition for the existence of a K\"ahler metric in the given Hermitian-symplectic Aeppli class.}

\vspace{1ex}

\section{Introduction}\label{section:Introduction} Let $X$ be an $n$-dimensional compact complex manifold. According to Sullivan [Sul76], Harvey-Lawson [HL83] and Streets-Tian [ST10, Definition 1.5] (where the name was coined), a Hermitian metric (namely a $C^\infty$ positive definite $(1,\,1)$-form) $\omega$ on $X$ is said to be {\bf Hermitian-symplectic} if $\omega$ is the component of bidegree $(1,1)$ of a real $C^\infty$ $d$-closed $2$-form $\widetilde{\omega}$ on $X$. Any $X$ admitting such a metric is called a {\bf Hermitian-symplectic manifold}. We will sometimes write H-S for Hermitian-symplectic.

These manifolds, which constitute a natural generalisation of compact K\"ahler manifolds, were given the following intrinsic characterisation by Sullivan.

\begin{The}\label{The:Sullivan_H-S} ([Sul76, Theorem III.2 and Remark III.11]) A compact complex manifold $X$ is Hermitian-symplectic if and only if $X$ carries no non-zero current $T$ of bidegree $(n-1,\,n-1)$ such that $T\geq 0$ and $T$ is $d$-exact.
\end{The}

Nevertheless, Hermitian-symplectic manifolds remain poorly understood. As they lie at the interface between symplectic and complex Hermitian geometries, they seem to warrant further probing. When $\mbox{dim}_\C X=2$, it can be shown (see e.g. [LZ09] or [ST10, Proposition 1.6] or Proposition \ref{Prop:S-T_surfaces} below) that $X$ is Hermitian-symplectic if and only if $X$ is K\"ahler. However, very little is known when $\mbox{dim}_\C X\geq 3$. This prompted Streets and Tian to ask the following

\begin{Question}\label{Question:S-T} ([ST10, Question 1.7]) Do there exist non-K\"ahler Hermitian-symplectic complex manifolds $X$ with $\mbox{dim}_\C X\geq 3$?

\end{Question}

While the general case of this question remains open, it has been answered negatively for a handful of special classes of manifolds, including all nilmanifolds endowed with an invariant complex structure by Enrietti, Fino and Vezzoni in [EFV12] and all twistor spaces by Verbitsky in [Ver14].

\vspace{2ex}

The Streets-Tian question is complementary to Donaldson's earlier

\begin{Question}\label{Question:Don} ([Don06, Question 2]) If $J$ is an almost-complex structure on a compact $4$-manifold which is tamed by a symplectic form, is there a symplectic form compatible with $J$?

\end{Question}

Indeed, when the almost-complex structure $J$ is integrable, a symplectic form $\widetilde\omega$ is a taming form for $J$ if and only if the $(1,\,1)$-component $\omega$ of $\widetilde\omega$ is a Hermitian-symplectic metric (i.e. positive definite). While $J$ is assumed integrable in Question \ref{Question:S-T}, the dimension of the underlying manifold is allowed to be arbitrary. Meanwhile, Question \ref{Question:Don}, that has come to be known in the literature as {\it Donaldson's tamed-to-compatible conjecture}, is peculiar to four real dimensions but $J$ need not be integrable. Thus, the only known case so far lies at the intersection of Questions \ref{Question:Don} and \ref{Question:S-T}.

\subsection{A new energy functional}\label{subsection:Introd_functional}

In this work, we investigate Question \ref{Question:S-T} by introducing a functional $F$ on the open convex subset ${\cal S}_{\{\omega_0\}}\subset\{\omega_0\}_A\cap C^\infty_{1,\,1}(X,\,\R)$ of all the Hermitian-symplectic metrics $\omega$ lying in the Aeppli cohomology class $\{\omega_0\}_A\in H^{1,\,1}_A(X,\,\R)$ of an arbitrary Hermitian-symplectic metric $\omega_0$.

Specifically, with every Hermitian-symplectic metric $\omega$ on $X$, we associate a unique differential form $\rho^{2,\,0}_\omega\in C^\infty_{2,\,0}(X,\,\C)$ that we call the {\bf $(2,\,0)$-torsion form} of $\omega$. (See Lemma and Definition \ref{Lem-Def:minimal_rho}.) We call its conjugate $\rho^{0,\,2}_\omega\in C^\infty_{0,\,2}(X,\,\C)$ the {\bf $(0,\,2)$-torsion form} of $\omega$. We then define our functional $F : {\cal S}_{\{\omega_0\}} \to [0,\,+\infty)$ (see Definition \ref{Def:F_energy-functional_H-S}) as the {\it squared $L^2$-norm} of $\rho^{2,\,0}_\omega$: \begin{equation}\label{eqn:Introd_F_energy-functional_H-S} F(\omega) = \int\limits_X|\rho^{2,\,0}_\omega|^2_\omega\,dV_\omega = ||\rho^{2,\,0}_\omega||^2_\omega.\end{equation}

  By construction, $F\geq 0$ and $F(\omega)=0$ if and only if the metric $\omega$ is K\"ahler. In the remaining part of $\S.$\ref{subsection:H-S_n=3_critical-points}, we go on to compute the {\it first variation} of $F$ in arbitrary dimension $n$ (see (ii) of Proposition \ref{Prop:F-tilde_properties}) and then reach the following conclusion in dimension $3$. (See also Corollary \ref{Cor:critical-points_energy_n=3}.)

  \begin{The}\label{The:Introd_critical-points_energy_n=3} Let $X$ be a $3$-dimensional compact Hermitian-symplectic manifold.

    For any H-S metric $\omega_0$ on $X$, the {\bf critical points} of $F : {\cal S}_{\{\omega_0\}} \to [0,\,+\infty)$ are the {\bf K\"ahler metrics} (if any) lying in the Aeppli cohomology class $\{\omega_0\}_A$.

  \end{The}

  In particular, the only possible critical points for $F$ are minima. Thus, the Streets-Tian Question \ref{Question:S-T} on $3$-dimensional compact complex manifolds $X$ is reduced to the existence of {\it critical points}, or equivalently {\it minimisers}, for the functional $F$.

\subsection{Generalised volume of Hermitian-symplectic Aeppli classes}\label{subsection:Introd_gen-volume}

Another consequence of Proposition \ref{Prop:F-tilde_properties} in the case of threefolds is that {\it minimising} our functional $F$ is equivalent to {\it maximising} the volume $\mbox{Vol}_\omega(X)$ of the Hermitian-symplectic metrics $\omega$ lying in a given Aeppli cohomology class. Specifically, we obtain the following result (see Proposition \ref{Prop:F-tilde_properties} and Definition \ref{Def:A-invariant}) that gives rise to a {\it volume-like invariant} for Aeppli cohomology classes of Hermitian-symplectic metrics.

\begin{The-Def}\label{The-Def:Introd_A-invariant} Let $X$ be a $3$-dimensional compact Hermitian-symplectic manifold.

  For any Hermitian-symplectic metric $\omega$ on $X$, the quantity \begin{equation}\label{eqn:Introd_A-invariant}A=A_{\{\omega\}_A}:= F(\omega) + \mbox{Vol}_\omega(X)>0\end{equation} is independent of the choice of metric $\omega$ in its Aeppli cohomology class $\{\omega\}_A$, where $\mbox{Vol}_\omega(X):=\int_X\omega^3/3!$.

    The invariant $A$ is called the {\bf generalised volume} of the Hermitian-symplectic Aeppli class $\{\omega\}_A$.

\end{The-Def}

Recall that, when $\omega$ is a K\"ahler metric (provided it exists) on an $n$-dimensional compact complex manifold $X$, the quantity $\mbox{Vol}_\omega(X)=\int_X\omega^n/n!$ depends only on the Bott-Chern class of $\omega$ and is standardly called the {\it volume} of the K\"ahler class $\{\omega\}_{BC}$ and denoted by $\mbox{Vol}(\{\omega\}_{BC})$. However, when $\omega$ is not K\"ahler, a major source of difficulty, for example in solving Monge-Amp\`ere equations, stems from $\mbox{Vol}_{\omega + i\partial\bar\partial\varphi}(X)$ depending on $i\partial\bar\partial\varphi$ when the real-valued smooth function $\varphi$ on $X$ varies such that $\omega + i\partial\bar\partial\varphi>0$. Thus, $\mbox{Vol}(\{\omega\}_{BC})$ is meaningless for non-K\"ahler classes, but it can be replaced by our generalised volume $A_{\{\omega\}_A}$, which coincides with the standard volume $\mbox{Vol}(\{\omega\}_{BC})$ when the class $\{\omega\}_{BC}$ is K\"ahler. See $\S.$\ref{subsection:vol-M-A_eq_H-S} for details on a new volume form and a natural Monge-Amp\`ere-type equation that we propose in the Hermitian-symplectic case.

\vspace{3ex}

A related problem appears in the study of compact {\bf SKT manifolds} i.e. compact complex manifolds admitting a Hermitian metric $\omega$ such that $\partial\bar{\partial}\omega=0$. It is well known that this class of manifolds is strictly larger than the K\"ahler class. For example, by [Gau77a], every compact complex surface is SKT.

It is of interest to investigate under what assumptions SKT manifolds are K\"ahler. In particular, the following special case of the Streets-Tian Question \ref{Question:S-T} is still open.

\begin{Question}\label{Question:SKT_ddbar} Let $X$ be a compact SKT manifold. Assume additionally that $X$ is a $\partial\bar{\partial}$-manifold. Does $X$ admit a K\"ahler metric?
\end{Question}

It is immediate to see that, on a $\partial\bar{\partial}$-manifold, a metric $\omega$ is Hermitian-symplectic if and only if $\omega$ is SKT. Our observation in $\S.$\ref{subsection:SKT_ddbar} is that a minor adaptation of the functional $F$ once again reduces the above problem to the existence of critical points.

\subsection{Search for critical points}\label{subsection:Introd_search-crit}

In complex dimension $3$, we propose a family of Monge-Amp\`ere-type equations and relate the maximisation of a numerical constant appearing therein to the local minimisers of the functional $F$.

\vspace{2ex}

Ideally, if a Monge-Amp\`ere-type equation with solutions in a given H-S {\bf Aeppli class} could be solved, its solutions would be K\"ahler metrics. Specifically, in $\S.$\ref{section:M-A_eq} we get the following result.

\begin{Prop}\label{Prop:Introd_consequence_MA-eq} Let $X$ be a compact complex Hermitian-symplectic manifold with $\mbox{dim}_\C X=3$. Fix an arbitrary Hermitian metric $\gamma$ and an H-S metric $\omega$ on $X$. Let $A = A_{\{\omega\}_A}>0$ be the generalised volume of the class $\{\omega\}_A$ and let $c = c_{\omega,\,\gamma}>0$ be the constant defined by the requirement \begin{equation*}\frac{(\int\limits_X\omega\wedge\gamma^2/2!)^3}{(\int\limits_X\gamma^3/3!)^2} = \frac{6A}{c}.\end{equation*}  

    If there exists a solution $\eta\in C^\infty_{1,\,0}(X,\,\C)$ of the  Monge-Amp\`ere-type equation \begin{equation*} (\omega + \partial\bar\eta + \bar\partial\eta)^3 = c\,(\Lambda_\gamma\omega)^3\,\,\frac{\gamma^3}{3!} \hspace{6ex} (\star)\end{equation*} such that $\omega_\eta:=\omega + \partial\bar\eta + \bar\partial\eta>0$, then $\omega_\eta$ is a {\bf K\"ahler metric} lying in the Aeppli cohomology class $\{\omega\}_A$ of $\omega$.

  \end{Prop}

\vspace{2ex}

However, equation $(\star)$ is heavily underdetermined and it is hard to see how one could go about solving it. For this reason, we replace it in $\S.$\ref{section:stratification} by a family of Monge-Amp\`ere-type equations of the familiar kind after we have stratified the given Hermitian-symplectic Aeppli class $\{\omega\}_A$ by {\it Bott-Chern subclasses} (or {\it strata}) in the following way.

We consider a {\it partition} of ${\cal S}_{[\omega]}$ of the shape ${\cal S}_{[\omega]}=\cup_{j\in J}{\cal D}_{[\omega_j]}$, where $(\omega_j)_{j\in J}$ is a family of H-S metrics in $\{\omega\}_A$ and $${\cal D}_{[\omega_j]}:=\{\omega'>0\,\mid\,\omega' - \omega_j\in\mbox{Im}\,(\partial\bar\partial)\}, \hspace{3ex} j\in J.$$

For each $j\in J$, we choose an arbitrary Hermitian metric $\gamma_j$ on $X$ such that $\Lambda_{\gamma_j}\omega_j=1$ at every point of $X$. Then, on each Bott-Chern stratum ${\cal D}_{[\omega_j]}\subset{\cal S}_{[\omega]}$, the Tosatti-Weinkove result in [TW10, Corollary 1] ensures the existence of a {\it unique} constant $b_j>0$ such that the equation \begin{equation*}\frac{(\omega_j + i\partial\bar\partial\varphi)^3}{3!} = b_jA\,\frac{dV_{\gamma_j}}{\int\limits_XdV_{\gamma_j}}\hspace{6ex} (\star\star\star_j),\end{equation*} subject to the extra condition $\omega_j + i\partial\bar\partial\varphi>0$, is {\it solvable}, where $A>0$ is the {\it generalised volume} of $\{\omega\}_A$ introduced in Theorem and Definition \ref{The-Def:Introd_A-invariant}. (Hence, $A$ is independent of $j$.)

In this way, we associate a constant $b_j\in(0,\,1]$ with every Bott-Chern stratum ${\cal D}_{[\omega_j]}$ of ${\cal S}_{[\omega]}$. The problem of {\it minimising} the functional $F$ in ${\cal S}_{[\omega]}$ (equivalently, {\it maximising} ${\cal S}_{[\omega]}\ni\omega'\mapsto\mbox{Vol}_{\omega'}(X)$) becomes equivalent to proving that the value $1$ is attained by one of the mysterious constants $b_j$. 

\begin{Prop}\label{Prop:Introd_b_j1} If there exists $j\in J$ such that $b_j=1$, the solution $\omega_j + i\partial\bar\partial\varphi_j$ of equation $(\star\star\star_j)$ is a {\bf K\"ahler metric} lying in the Bott-Chern subclass ${\cal D}_{[\omega_j]}$, hence also in the Aeppli class $\{\omega\}_A$.

\end{Prop}

We go on to observe in Lemma \ref{Lem:BC-subclass_volume-Gauduchon} that if a Bott-Chern stratum of an H-S Aeppli class contains a {\it Gauduchon} metric, then all the metrics $\omega$ on that stratum are Gauduchon and have the same volume $\mbox{Vol}_\omega(X)$. On the other hand, the restriction of the volume function $\omega\mapsto\mbox{Vol}_\omega(X)$ to a non-Gauduchon stratum cannot achieve any local extremum, thanks to Lemma \ref{Lem:BC-subclass_volume-nonGauduchon}. Thus, we have a good understanding of the behaviour of the volume in the {\it horizontal} directions (i.e. those of the Bott-Chern strata). The variation in the {\it vertical} directions remains mysterious for now.

\subsection{Obstruction to the existence of a K\"ahler metric in a given Hermitian-symplectic Aeppli class}\label{subsection:Introd_obstruction}

While the Streets-Tian Question \ref{Question:S-T} asks whether a K\"ahler metric exists on every Hermitian-symplectic manifold $X$, if our functional $F : {\cal S}_{\{\omega_0\}} \to [0,\,+\infty)$ admits critical points for {\it any} Hermitian-symplectic metric $\omega_0$ on $X$, much more will be true: there will exist a K\"ahler metric in the Aeppli cohomology class of {\it every} Hermitian-symplectic metric on $X$.  

  However, in $\S.$\ref{section:obs-est} we identify an obstruction to the existence of a  K\"ahler metric that is Aeppli cohomologous to a given Hermitian-symplectic metric.

  In fact, we first show (see Lemma and Definition \ref{Def:E_2_H-S}) that the $(0,\,2)$-torsion form $\rho^{0,\,2}_\omega$ of any H-S metric $\omega$ on an $n$-dimensional compact complex manifold $X$ is {\bf $E_2$-closed} in the sense that it defines an $E_2$-cohomology class $$\{\rho^{0,\,2}_\omega\}_{E_2}\in E_2^{0,\,2}(X)$$ on the second page of the {\it Fr\"olicher spectral sequence} of $X$. Moreover, $\{\rho^{0,\,2}_\omega\}_{E_2}$ depends only on the Aeppli class $\{\omega\}_A$. We call it the {\bf $E_2$-torsion class} of the Hermitian-symplectic Aeppli class $\{\omega\}_A$ and prove the following fact (see Corollary \ref{Cor:necessary-cond_K}).

  \begin{Prop}\label{Prop:Introd_} Let $X$ be a $3$-dimensional compact complex manifold supposed to carry a Hermitian-symplectic metric $\omega$.

    The {\bf vanishing} of the $E_2$-torsion class $\{\rho^{0,\,2}_\omega\}_{E_2}\in E_2^{0,\,2}(X)$ is a necessary condition for the Aeppli class $\{\omega\}_A$ to contain a K\"ahler metric.

  \end{Prop}

  This throws up the natural question of whether there exist compact $3$-dimensional Hermitian-symplectic manifolds on which all or some $E_2$-torsion classes are non-vanishing.

  \subsection{Cohomological interpretations of the generalised volume}\label{subsection:Introd_coh-interpretations}

In $\S.$\ref{subsection:cohom_A} and $\S.$\ref{subsection:minimal-completion}, we give two cohomological interpretations of our generalised volume invariant $A$.

\vspace{2ex}

The one in $\S.$\ref{subsection:minimal-completion} is valid on any $3$-dimensional compact Hermitian-symplectic manifold. We first observe that the real $d$-closed $2$-form $$\widetilde\omega = \rho_\omega^{2,\,0} + \omega + \rho_\omega^{0,\,2},$$ that we call the {\bf minimal completion} of the given H-S metric $\omega$, lies in a De Rham cohomology class that depends only on the Aeppli cohomology class of $\omega$. This follows from the study of our functional $F$, specifically from Corollary \ref{Cor:energy_M-A_mass}. We then observe (see (a) of Proposition \ref{Prop:same-DR-class}) that the {\it generalised volume} is a kind of {\it volume of the minimal completion}:  \begin{equation}\label{eqn:Introd_min-comp_integral}A = A_{\{\omega\}_A} = \int\limits_X\frac{\widetilde\omega^3}{3!} = \frac{1}{6}\,\{\widetilde\omega\}_{DR}^3,\end{equation} so it only depends on the De Rham cohomology class of the minimal completion $\widetilde\omega$.

\vspace{2ex}

The cohomological interpretation given in $\S.$\ref{subsection:cohom_A} is only valid on $3$-dimensional compact Hermitian-symplectic manifolds that lie in the new class of {\bf page-$1$-$\partial\bar\partial$-manifolds} that were very recently introduced in [PSU20a]. This class is strictly larger than the one of $\partial\bar\partial$-manifolds, so we get a link with Question \ref{Question:SKT_ddbar}.

As observed in [PSU20b, Proposition 6.2], on a $3$-dimensional manifold $X$ an H-S metric $\omega$ induces an {\it $E_2$-Aeppli cohomology} class $\{\omega\}_{E_2,\,A}\in E_{2,\,A}^{1,\,1}(X)$. Together with the {\it $E_r$-Bott-Chern cohomologies}, the {\it $E_r$-Aeppli cohomologies} have been recently introduced in [PSU20b, Definition 3.4] for every integer $r\geq 2$. They coincide with the standard Bott-Chern and Aeppli cohomologies when $r=1$. 

On the other hand, using results from [PSU20b], we show in Corollary \ref{Cor:E_2BC_lifts} that on a $3$-dimensional page-$1$-$\partial\bar\partial$-manifold $X$, an {\it $E_2$-Bott-Chern class} $\mathfrak{c}_\omega\in E_{2,\,BC}^{2,\,2}(X)$ can be canonically associated with the $E_2$-Aeppli class $\{\omega\}_{E_2,\,A}$ of any Hermitian-symplectic metric $\omega$ on $X$. Finally, using the duality between $E_{2,\,BC}^{n-1,\,n-1}(X)$ and $E_{2,\,A}^{1,\,1}(X)$ proved in [PSU20b, Theorem 3.11] for every compact complex manifold $X$ of any dimension $n$, we get the following cohomological interpretation (see Theorem \ref{The:cohom_A}) of the generalised volume as a multiple of the {\it intersection number} between the cohomology classes $\mathfrak{c}_\omega$ and $\{\omega\}_{E_2,\,A}$: \begin{equation}\label{eqn:Introd_cohom_A} A= A_{\{\omega\}_A} = \frac{1}{6}\,\mathfrak{c}_\omega.\{\omega\}_{E_2,\,A}.\end{equation}

Corollary \ref{Cor:energy_M-A_mass} in the study of our functional $F$ is again used in a key way to obtain this result.

\vspace{3ex}

\noindent {\bf Acknowledgments.} This work started in the spring of 2016 when the second-named author was visiting the Jagiellonian University in Krak\'ow at the invitation of S\l{}awomir Kolodziej, to whom he is very grateful for the hospitality offered under the NCN grant 2013/08/A/ST1/00312 and the very stimulating discussions on various topics. The first-named author was partially supported by  the  National Science Centre, Poland grant no 2017/27/B/ST1/01145.

 \section{Preliminaries}\label{section:preliminaries} In this section, we present a mixture of well-known and new results that will come in handy.

\subsection{Background}\label{subsection:background} Let $\omega$ be a Hermitian metric on a compact complex manifold $X$ with $\mbox{dim}_\C X=n$. To outline the context of our study, we start by recalling the definitions of six classes of special metrics and the known implications among them: \\

\vspace{3ex}

\noindent$\begin{array}{lllll}d\omega=0 & \Longrightarrow & \exists\,\, \rho^{0,\,2}\in C^{\infty}_{0,\,2}(X,\,\C)\,\, \mbox{s.t.} & \Longrightarrow & \partial\bar\partial\omega=0 \\
 &  & d(\overline{\rho^{0,\,2}}+\omega+\rho^{0,\,2})=0 & & \\
 (\omega\,\,\mbox{is K\"ahler}) &   & (\omega\,\,\mbox{is Hermitian-symplectic}) &  & (\omega\,\,\mbox{is SKT}) \\
\rotatebox{-90}{$\implies$} &  & & & \hspace{33ex} (P) \\
d\omega^{n-1}=0 & \Longrightarrow & \exists\,\, \Omega^{n-2,\,n}\in C^{\infty}_{n-2,\,n}(X,\,\C)\,\, \mbox{s.t.}  & \Longrightarrow & \partial\bar\partial\omega^{n-1}=0 \\
&  &  d(\overline{\Omega^{n-2,\,n}}+\omega^{n-1}+\Omega^{n-2,\,n})=0 & &   \\
(\omega\,\,\mbox{is balanced}) &   & (\omega\,\,\mbox{is strongly Gauduchon (sG)}) &  & (\omega\,\,\mbox{is Gauduchon}).\end{array}$

\vspace{5ex}

\noindent The manifold $X$ is called {\it K\"ahler}, {\it Hermitian-symplectic} (H-S), {\it SKT}, {\it balanced}, {\it strongly Gauduchon} (sG) if it carries a Hermitian metric $\omega$ of the corresponding type.  Meanwhile, Gauduchon metrics always exist on any $X$ by [Gau77a]. Balanced metrics were introduced in [Gau77b] under the name {\it semi-K\"ahler} and then discussed again in [Mic83], while strongly Gauduchon (sG) metrics were introduced in [Pop13] by requiring $\partial\omega^{n-1}\in\mbox{Im}\,\bar\partial$, a definition that was then proved in [Pop13, Proposition 4.2] to be equivalent to the description on the second line in the above picture (P). In particular, the notion of H-S metric is the analogue in bidegree $(1,\,1)$ of the notion of sG metric. 

These special metrics define Dolbeault, Bott-Chern or Aeppli cohomology classes, according to the case. Recall the by now standard definitions of these cohomologies: $$H^{\bullet,\,\bullet}_{\bar\partial}(X,\,\C)=\ker\bar\partial/\mbox{Im}\,\bar\partial, \hspace{1ex}  H^{\bullet,\,\bullet}_{BC}(X,\,\C)=\ker\partial\cap\ker\bar\partial/\mbox{Im}\,(\partial\bar\partial), \hspace{1ex} H^{\bullet,\,\bullet}_A(X,\,\C)=\ker(\partial\bar\partial)/(\mbox{Im}\,\partial + \mbox{Im}\,\bar\partial).$$

The metrical notions in bidegree $(n-1,\, n-1)$ listed on the second line in (P) offer more flexibility than their bidegree $(1,\, 1)$ analogues. Moreover, the cohomology classes they define eventually give information on those defined by the metrics on the first line thanks to the classical {\it Serre duality} for the Dolbeault cohomology and the analogous duality between the Bott-Chern and Aeppli cohomologies (see [Sch07] for the latter):

$$H^{p,\,q}_{\bar\partial}(X,\,\C)\times H^{n-p,\,n-q}_{\bar\partial}(X,\,\C)\longrightarrow\C, \hspace{3ex} (\{\alpha\}_{\bar\partial},\,(\{\beta\}_{\bar\partial})\longmapsto\int\limits_X\alpha\wedge\beta,$$

\noindent and

$$H^{p,\,q}_{BC}(X,\,\C)\times H^{n-p,\,n-q}_A(X,\,\C)\longrightarrow\C, \hspace{3ex} (\{\alpha\}_{BC},\,(\{\beta\}_A)\longmapsto\int\limits_X\alpha\wedge\beta.$$

As for other possible vertical implications in (P), besides the trivial ``$\omega$ {\it K\"ahler $\implies$ $\omega$ balanced}'', it is easy to see that there is no counterpart at the SKT/Gauduchon level:

\vspace{1ex}

\hspace{25ex} $\omega$ {\it SKT $\centernot\implies$ $\omega$ Gauduchon}.

\subsection{H-S and sG manifolds}\label{subsection:H-S_sG} We will now prove the implication ``{\it  Hermitian-symplectic $\implies$ strongly Gauduchon}'' at the level of manifolds. It was given a very similar proof as Lemma $1$ in $\S.2$ of [YZZ19], that we now recall for the reader's convenience. It goes some way in the direction of the Streets-Tian Question \ref{Question:S-T}. Note that this implication does not hold at the level of metrics $\omega$.

\begin{Prop}\label{Prop:H-S_sG} Every compact complex manifold $X$ that admits a {\bf Hermitian-symplectic} metric also admits a {\bf strongly Gauduchon (sG)} metric. 

\end{Prop}  
 
\noindent {\it Proof.} Let $n=\mbox{dim}_\C X$. As recalled in (P), a strongly Gauduchon (sG) structure on $X$ can be regarded as a real $C^\infty$ $d$-closed $(2n-2)$-form $\Omega$ on $X$ such that its $(n-1,\,n-1)$-component $\Omega^{n-1,\,n-1}$ is positive definite (see [Pop13]). This also uses the fact that, if $\Omega^{n-1,\,n-1}>0$, there exists a unique smooth positive definite $(1,\,1)$-form $\omega$ on $X$, called the $(n-1)$-st root $\Omega^{n-1,\,n-1}$, such that $\omega^{n-1} = \Omega^{n-1,\,n-1}$. (This fact, noticed in [Mic83], is well known and can be easily checked pointwise in appropriately chosen local coordinates.)

Now, suppose that an H-S structure $\widetilde\omega$ exists on $X$. This means that $\widetilde\omega = \rho^{2,\,0} + \omega + \rho^{0,\,2}$ is a real $C^\infty$ $d$-closed $2$-form on $X$ such that its $(1,\,1)$-component $\omega$ is positive definite. Thus, $d\widetilde\omega^{n-1}=0$ and \begin{eqnarray}\nonumber\widetilde\omega^{n-1} = [\omega + (\rho^{2,\,0} + \rho^{0,\,2})]^{n-1} = \sum\limits_{k=0}^{n-1}\sum\limits_{l=0}^k{n-1 \choose k}{k \choose l}\,(\rho^{2,\,0})^l\wedge(\rho^{0,\,2})^{k-l}\wedge\omega^{n-k-1}.\end{eqnarray}

\noindent In particular, the $(n-1,\,n-1)$-component of $\widetilde\omega^{n-1}$ is the sum of the terms for which $l=k-l$ in the above expression, i.e. $$\Omega^{n-1,\,n-1} = \omega^{n-1} + \sum\limits_{l=1}^{[\frac{n-1}{2}]}{n-1 \choose 2l}{2l \choose l}\,(\rho^{2,\,0})^l\wedge(\rho^{0,\,2})^l\wedge\omega^{n-2l-1}.$$

Thus, to prove the existence of an sG structure on $X$, it suffices to prove that the $(n-1,\,n-1)$-form $\Omega^{n-1,\,n-1}$ is positive definite. Its $(n-1)$-st root will then be an sG metric on $X$, by construction.

To show that $\Omega^{n-1,\,n-1}>0$, it suffices to check that the real form $(\rho^{2,\,0})^l\wedge(\rho^{0,\,2})^l\wedge\omega^{n-2l-1}$ is weakly (semi)-positive at every point of $X$. (Recall that $\rho^{0,\,2}$ is the conjugate of $\rho^{2,\,0}$.) To this end, note that the $(2l,2l)$-form $(\rho^{2,\,0})^l\wedge(\rho^{0,\,2})^l$ is weakly semi-positive as the wedge product of a $(2l,0)$-form and its conjugate (see [Dem97, Chapter III, Example 1.2]). Therefore, the $(n-1,\,n-1)$-form $(\rho^{2,\,0})^l\wedge(\rho^{0,\,2})^l\wedge\omega^{n-2l-1}$ is (semi)-positive since the product of a weakly (semi)-positive form and a strongly (semi)-positive form is weakly (semi)-positive and $\omega$ is strongly positive (see [Dem97, Chapter III, Proposition 1.11]). (Recall that in bidegrees $(1,\,1)$ and $(n-1,\,n-1)$, the notions of weak and strong positivity coincide.)\hfill $\Box$

\vspace{3ex}

In the case $n=2$, the notions of H-S and sG metrics coincide. Meanwhile, as explained in [Pop13, Observation 4.4], every strongly Gauduchon compact complex surface is K\"ahler. This answers the two-dimensional analogue of the Streets-Tian Question \ref{Question:S-T}, a fact that has been known for a while (cf. e.g. [LZ09] or [ST10, Proposition 1.6]).

\begin{Prop}\label{Prop:S-T_surfaces} Let $X$ be a compact complex surface. If $X$ carries a Hermitian-symplectic metric, then $X$ carries a K\"ahler metric.

\end{Prop}  

Note that Theorem 6.1 in [Lam99], according to which every compact complex surface whose first Betti number $b_1$ is odd carries a non-zero positive $d$-exact $(1,\,1)$-current, is given a proof based on the Hahn-Banach separation theorem and uses a duality argument. On the other hand, Observation 4.4 in [Pop13] uses in a key way Lamari's result and the fact that a compact complex manifold carries a strongly Gauduchon metric if and only there is no non-zero positive $d$-exact $(1,\,1)$-current on it (see [Pop13, Proposition 4.3]). This characterisation of sG manifolds given in [Pop13] uses again the Hahn-Banach theorem and a Harvey-Lawson-type duality argument harking back to Sullivan. In particular, we get no information in this way on the Aeppli cohomology class of the K\"ahler metric whose existence is given by the above Proposition \ref{Prop:S-T_surfaces}.

\vspace{3ex}

Independently, let us notice that the existence of Hermitian-symplectic metrics on a compact complex threefold implies a property that is well known to hold on compact K\"ahler manifolds (and even on $\partial\bar\partial$-manifolds and even on compact complex manifolds whose Fr\"olicher spectral sequence degenerates at $E_1$). Thus, the observation in Corollary \ref {Cor:holomorphic_1-forms} takes Hermitian-symplectic threefolds a little closer to K\"ahler ones.

In the proof of  Corollary \ref {Cor:holomorphic_1-forms} and thereafter, we will need the following standard formula (cf. e.g. [Voi02, Proposition 6.29, p. 150]) for the Hodge star operator $\star = \star_\omega$ of any Hermitian metric $\omega$ applied to {\it primitive} forms $v$ of arbitrary bidegree $(p, \, q)$: \begin{eqnarray}\label{eqn:prim-form-star-formula-gen}\star\, v = (-1)^{k(k+1)/2}\, i^{p-q}\, \frac{\omega^{n-p-q}\wedge v}{(n-p-q)!}, \hspace{2ex} \mbox{where}\,\, k:=p+q.\end{eqnarray} Recall that, for a given integer $k\in\{0,\dots , n\}$, a $k$-form $v$ on an $n$-dimensional complex manifold is said to be {\it primitive} (w.r.t. a fixed Hermitian metric $\omega$) if $\omega^{n-k+1}\wedge v =0$. This is known to be equivalent to $\Lambda_\omega v=0$, where $\Lambda_\omega$ is the adjoint of the Lefschetz operator $\omega\wedge\cdot$ w.r.t. the pointwise inner product $\langle\,\,,\,\,\rangle_\omega$ defined by $\omega$. In particular, all $k$-forms with $k\in\{0,\,1\}$ are primitive and so are all $(p,\,0)$-forms and all $(0,\,q)$-forms.

\vspace{3ex}

\begin{Cor}\label{Cor:holomorphic_1-forms} Let $X$ be a compact complex Hermitian-symplectic manifold with $\mbox{dim}_\C X=3$. 

 Then, every holomorphic $1$-form (i.e. every smooth $\bar\partial$-closed $(1,\,0)$-form)  on $X$ is $d$-closed.

\end{Cor}

\noindent {\it Proof.} Let $\omega$ be an H-S metric on $X$. Then, $\partial\omega\in\mbox{Im}\,\bar\partial$ and $\bar\partial\omega\in\mbox{Im}\,\partial$ (see e.g. (\ref{eqn:H-S_condition})). Choose any form $\rho^{2,\,0}\in C^\infty_{2,\,0}(X,\,\C)$ such that $\partial\omega = -\bar\partial\rho^{2,\,0}$. Hence, $\bar\partial\omega = -\partial\rho^{0,\,2}$, where $\rho^{0,\,2}:=\overline{\rho^{2,\,0}}$.

Now, let $\xi\in C^{\infty}_{1,\,0}(X,\,\C)$ such that $\bar\partial\xi=0$. We want to show that $\partial\xi=0$.

On the one hand, if $\star = \star_\omega$ is the Hodge star operator induced by $\omega$, the general formula (\ref{eqn:prim-form-star-formula-gen}) applied to the (necessarily primitive) $(0,\,2)$-form $\bar\partial\bar\xi$ yields: $\star(\bar\partial\bar\xi) = \bar\partial\bar\xi\wedge\omega$. Hence, \begin{equation}\label{eqn:prelim_d-closedness1}\partial\xi\wedge\bar\partial\bar\xi\wedge\omega = |\partial\xi|_\omega^2\,dV_\omega \geq 0\end{equation} at every point of $X$.

  Meanwhile, an immediate calculation and the use of the identities $\bar\partial\xi=0$ and $\partial\bar\xi=0$ show that \begin{eqnarray*}\partial\xi\wedge\bar\partial\bar\xi\wedge\omega & = & -\partial\bar\partial(\xi\wedge\bar\xi\wedge\omega) + \xi\wedge\bar\partial\bar\xi\wedge\partial\omega + \partial\xi\wedge\bar\xi\wedge\bar\partial\omega + \xi\wedge\bar\xi\wedge\partial\bar\partial\omega \\
    & = & -\partial\bar\partial(\xi\wedge\bar\xi\wedge\omega) - \xi\wedge\bar\partial\bar\xi\wedge\bar\partial\rho^{2,\,0} - \partial\xi\wedge\bar\xi\wedge\partial\rho^{0,\,2} \\
    & = & -\partial\bar\partial(\xi\wedge\bar\xi\wedge\omega) + \bar\partial(\xi\wedge\bar\partial\bar\xi\wedge\rho^{2,\,0}) + \partial(\partial\xi\wedge\bar\xi\wedge\rho^{0,\,2})\in\mbox{Im}\,\partial + \mbox{Im}\,\bar\partial,\end{eqnarray*} where for the second identity we also used the property $\partial\bar\partial\omega=0$ of the H-S metric $\omega$. Using Stokes's theorem, we infer: \begin{equation}\label{eqn:prelim_d-closedness2}\int\limits_X\partial\xi\wedge\bar\partial\bar\xi\wedge\omega =  0.\end{equation}

  Putting together (\ref{eqn:prelim_d-closedness1}) and (\ref{eqn:prelim_d-closedness2}), we get $\partial\xi=0$ on $X$ and we are done.  \hfill $\Box$

\subsection{Toolbox}\label{subsection:toolbox}

\hspace{2ex} (I)\, Let $\omega$ be an arbitrary Hermitian metric on an $n$-dimensional compact complex manifold $X$. We know from [KS60, $\S.6$] (see also [Sch07] or [Pop15]) that $\omega$ induces the following $4^{th}$-order elliptic differential operator $\Delta_{BC}: C^{\infty}_{r,\,s}(X,\,\C)\longrightarrow C^{\infty}_{r,\,s}(X,\,\C)$, called the {\it Bott-Chern Laplacian}, in every bidegree $(r,\,s)$: \begin{equation}\label{eqn:BC-Laplacian}\Delta_{BC} : = \partial^{\star}\partial + \bar\partial^{\star}\bar\partial + (\partial\bar\partial)^{\star}(\partial\bar\partial) + (\partial\bar\partial)(\partial\bar\partial)^{\star} + (\partial^{\star}\bar\partial)^{\star}(\partial^{\star}\bar\partial) + (\partial^{\star}\bar\partial)(\partial^{\star}\bar\partial)^{\star}.\end{equation}

\noindent From the ellipticity and self-adjointness of $\Delta_{BC}$, coupled with the compactness of $X$, we get the following $L^2_\omega$-orthogonal $3$-space decomposition: \begin{equation}\label{eqn:BC-3sp-decomp}C^{\infty}_{r, \, s}(X, \C)=\ker\Delta_{BC} \oplus \mbox{Im}\,\partial\bar\partial \oplus (\mbox{Im}\,\partial^{\star} + \mbox{Im}\,\bar\partial^{\star})\end{equation}

\noindent in which $\ker\partial\cap\ker\bar\partial = \ker\Delta_{BC} \oplus \mbox{Im}\,\partial\bar\partial$.

\vspace{3ex}

The following Neumann-type formula for the minimal $L^2_\omega$-norm solution of a $\bar\partial$-equation with an extra constraint seems to be new. It will be used in $\S.$\ref{section:obs-est}.

\begin{Lem}\label{Lem:Neumann_del-0} Let $(X,\,\omega)$ be a compact Hermitian manifold. For every $p,q=0,\dots , n=\mbox{dim}_\C X$ and every form $v\in C^\infty_{p,\,q}(X,\,\C)$, consider the following $\bar\partial$-equation problem: \begin{equation}\label{eqn:Neumann_del-0}\bar\partial u = v  \hspace{3ex} \mbox{subject to the condition} \hspace{2ex} \partial u = 0.\end{equation}

\noindent If problem (\ref{eqn:Neumann_del-0}) is solvable for $u$, the solution of minimal $L^2_\omega$-norm is given by the Neumann-type formula: \begin{equation}\label{eqn:Neumann_formula_del-0}u =\Delta_{BC}^{-1}[\bar\partial^\star v + \bar\partial^\star\partial\partial^\star v].\end{equation}

\end{Lem}

\noindent {\it Proof.} The solution $u$ of problem (\ref{eqn:Neumann_del-0}) is unique up to $\ker\partial\cap\ker\bar\partial= \ker\Delta_{BC} \oplus \mbox{Im}\,\partial\bar\partial$. Thanks to (\ref{eqn:BC-3sp-decomp}), the minimal $L^2_\omega$-norm solution of problem (\ref{eqn:Neumann_del-0}) is uniquely determined by the condition $u\in\mbox{Im}\,\partial^{\star} + \mbox{Im}\,\bar\partial^{\star}$. In other words, there exist forms $\xi$ and $\eta$ such that 

$$u = \partial^\star\xi + \bar\partial^\star\eta, \hspace{3ex} \mbox{hence} \hspace{3ex} \partial^{\star}u = -\bar\partial^{\star}\partial^{\star}\eta, \hspace{2ex} \bar\partial^{\star}u = -\partial^{\star}\bar\partial^{\star}\xi \hspace{2ex} \mbox{and} \hspace{2ex} (\partial\bar\partial)^{\star}u = 0.$$

\noindent Applying $\Delta_{BC}$, we get

$$\Delta_{BC} u = \bar\partial^{\star}(\bar\partial u) + \bar\partial^{\star}\partial\partial^\star(\bar\partial u),$$

\noindent since the first, third (after writing $\partial\bar\partial = -\bar\partial\partial$) and sixth (after writing $(\partial^{\star}\bar\partial)^{\star} = \bar\partial^\star\partial$) terms in $\Delta_{BC}$ end with $\partial$ and $\partial u=0$, while the fourth term in $\Delta_{BC}$ ends with $(\partial\bar\partial)^\star$ and $(\partial\bar\partial)^\star u=0$.

Now, the restriction of $\Delta_{BC}$ to the orthogonal complement of $\ker\Delta_{BC}$ is an isomorphism onto this same orthogonal complement, so using the inverse $\Delta_{BC}^{-1}$ of this restriction ($=$ the Green operator of $\Delta_{BC}$), we get

$$u=\Delta_{BC}^{-1}[\partial^{\star}(\partial u) + \partial^{\star}\bar\partial\bar\partial^\star(\partial u)],$$

\noindent since both $u$ and $\partial^{\star}(\partial u) + \partial^{\star}\bar\partial\bar\partial^\star(\partial u)$ lie in $(\ker\Delta_{BC})^\perp$. 

Since $\partial u =v$, the last formula for $u$ is precisely (\ref{eqn:Neumann_formula_del-0}). \hfill $\Box$

\vspace{3ex}

(II)\, We now remind the reader of the following notion.

\begin{Def}\label{Def:dd-bar-lemma} A compact complex manifold $X$ is said to be a $\partial\bar\partial$-{\bf manifold} if for any $d$-closed {\it pure-type} form $u$ on $X$, the following exactness properties are equivalent: \\

\hspace{10ex} $u$ is $d$-exact $\Longleftrightarrow$ $u$ is $\partial$-exact $\Longleftrightarrow$ $u$ is $\bar\partial$-exact $\Longleftrightarrow$ $u$ is $\partial\bar\partial$-exact.

\end{Def}

Recall that $\partial\bar\partial$-manifolds are precisely the compact complex manifolds that have the {\it canonical Hodge Decomposition property} in the sense that every Dolbeault cohomology class $\{\alpha\}_{\bar\partial}\in H^{p,\,q}_{\bar\partial}(X,\,\C)$ of any bidegree $(p,\,q)$ can be represented by a {\it $d$-closed} form and the identity map induces, via $d$-closed representatives, an {\it isomorphism} $$H^k_{DR}(X,\,\C)\simeq\bigoplus\limits_{p+q=k}H^{p,\,q}_{\bar\partial}(X,\,\C)$$ in every degree $k$.

\vspace{3ex}

(III)\, Let us also mention the following observation that will play a key part in this work. It was first noticed in [IP13] and in some of the references therein as a consequence of more general results. A quick proof, made even shorter below, appeared in [Pop15, Proposition 1.1]. 

\begin{Prop}\label{Prop:SKT+bal} If a Hermitian metric $\omega$ on a compact complex manifold $X$ is both {\bf SKT} and {\bf balanced}, then $\omega$ is {\bf K\"ahler}.

\end{Prop}

\noindent {\it Proof.} The {\it SKT} assumption on $\omega$ translates to any of the following equivalent properties: \begin{eqnarray}\label{eqn:pluriclosed-equiv}\partial\bar\partial\omega=0 \Longleftrightarrow \partial\omega\in\ker\bar\partial \Longleftrightarrow \star(\partial\omega)\in\ker\partial^{\star},\end{eqnarray}

\noindent where the last equivalence follows from the standard formula $\partial^{\star}=-\star\bar\partial\star$ involving the Hodge-star isomorphism $\star=\star_{\omega}:\Lambda^{p,\,q}T^{\star}X\rightarrow \Lambda^{n-q,\,n-p}T^{\star}X$ defined by $\omega$ for arbitrary $p,q=0,\dots , n$.

 Meanwhile, the {\it balanced} assumption on $\omega$ translates to any of the following equivalent properties: \begin{eqnarray*}\label{eqn:bal-equiv}d\omega^{n-1}=0 \Longleftrightarrow \partial\omega^{n-1}=0 \Longleftrightarrow \omega^{n-2}\wedge\partial\omega = 0   \Longleftrightarrow \partial\omega \,\,\mbox{is primitive}.\end{eqnarray*}

\noindent Moreover, since $\partial\omega$ is primitive when $\omega$ is balanced, the general formula (\ref{eqn:prim-form-star-formula-gen}) yields: \begin{eqnarray}\label{eqn:consequence_balanced}\star(\partial\omega) = i\,\frac{\omega^{n-3}}{(n-3)!}\wedge\partial\omega = \frac{i}{(n-2)!}\,\partial\omega^{n-2}\in\mbox{Im}\,\partial.\end{eqnarray}

\vspace{1ex}

Thus, if $\omega$ is both SKT and balanced, we get from (\ref{eqn:pluriclosed-equiv}) and (\ref{eqn:consequence_balanced}) that \begin{eqnarray*}\star(\partial\omega)\in\ker\partial^{\star}\cap\mbox{Im}\,\partial = \{0\},\end{eqnarray*} where the last identity follows from the subspaces $\ker\partial^{\star}$ and $\mbox{Im}\,\partial$ of $C^\infty_{n-1,\,n-2}(X,\,\C)$ being $L^2_\omega$-orthogonal. We infer that $\partial\omega=0$, i.e. $\omega$ is K\"ahler.   \hfill  $\Box$

\section{The energy functional}\label{section:alternative_functional}

We will define and discuss our new energy functional in the general Hermitian-symplectic setting in $\S.$\ref{subsection:H-S_n=3_critical-points}. We will then discuss a variant of it in the special case of SKT $\partial\bar\partial$-manifolds in $\S.$\ref{subsection:SKT_ddbar}. Additional discussions are included in $\S.$\ref{subsection:variation_torsion} and $\S.$\ref{subsection:vol-M-A_eq_H-S}.

\subsection{Case of H-S metrics on  compact complex manifolds}\label{subsection:H-S_n=3_critical-points}

Let $X$ be a compact complex manifold with $\mbox{dim}_{\C}X=n$ such that $X$ admits {\it Hermitian-symplectic} metrics. Recall that these are $C^{\infty}$ positive definite $(1,\,1)$-forms $\omega>0$ for which there exists $\rho^{2,\,0}\in C^{\infty}_{2,\,0}(X,\,\C)$ such that 

\begin{equation}\label{eqn:H-S_condition}d(\rho^{2,\,0} + \omega + \rho^{0,\,2}) = 0,\end{equation} where $\rho^{0,\,2}:=\overline{\rho^{2,\,0}}$. Alternatively, we say that $\widetilde{\omega}:=\rho^{2,\,0} + \omega + \rho^{0,\,2}$ is a Hermitian-symplectic $2$-form.

\begin{Lem-Def}\label{Lem-Def:minimal_rho} For every Hermitian-symplectic metric $\omega$ on $X$, there exists a unique smooth $(2,\,0)$-form $\rho_\omega^{2,\,0}$ on $X$ such that

 \begin{equation}\label{eqn:H-S_condition_bis}(i)\,\, \partial\rho_\omega^{2,\,0} = 0  \hspace{3ex} \mbox{and} \hspace{3ex}  (ii)\,\, \bar\partial\rho_\omega^{2,\,0} = -\partial\omega \hspace{3ex} \mbox{and} \hspace{3ex} (iii)\,\,\rho_\omega^{2,\,0}\in \mbox{Im}\,\partial^{\star}_\omega + \mbox{Im}\,\bar\partial^{\star}_\omega.\end{equation}

\noindent Moreover, property $(iii)$ ensures that $\rho_\omega^{2,\,0}$ has {\bf minimal $L^2_\omega$ norm} among all the $(2,\,0)$-forms satisfying properties $(i)$ and $(ii)$.

 We call $\rho_\omega^{2,\,0}$ the {\bf $(2,\,0)$-torsion form} and its conjugate $\rho_\omega^{0,\,2}$ the {\bf $(0,\,2)$-torsion form} of the Hermitian-symplectic metric $\omega$. One has the explicit {\bf Neumann-type formula}:

\begin{equation}\label{eqn:Neumann_torsion_H-S}\rho_\omega^{2,\,0} = -\Delta_{BC}^{-1}[\bar\partial^\star\partial\omega + \bar\partial^\star\partial\partial^\star\partial\omega],\end{equation}

\noindent where $\Delta_{BC}^{-1}$ is the Green operator of the Bott-Chern Laplacian $\Delta_{BC}$ induced by $\omega$, while $\partial^\star=\partial^\star_\omega$ and $\bar\partial^\star=\bar\partial^\star_\omega$ are the formal adjoints of $\partial$, resp. $\bar\partial$, w.r.t. the $L^2$ inner product defined by $\omega$.

\end{Lem-Def}

\noindent {\it Proof.} Condition (\ref{eqn:H-S_condition}) is equivalent to the vanishing of each of the components of pure types $(3,\,0)$, $(2,\,1)$, $(1,\,2)$ and $(0,\,3)$ of the real $3$-form $d(\rho^{2,\,0} + \omega + \rho^{0,\,2})$. Since the $(3,\,0)$- and $(2,\,1)$-components are the conjugates of the $(0,\,3)$- and resp. $(1,\,2)$-components, these vanishings are equivalent to conditions $(i)$ and $(ii)$ of (\ref{eqn:H-S_condition_bis}) being satisfied by $\rho^{2,\,0}$ in place of $\rho^{2,\,0}_\omega$. 

Now, the forms $\rho^{2,\,0}$ satisfying equations $(i)$ and $(ii)$ of (\ref{eqn:H-S_condition_bis}) are unique modulo $\ker\partial\cap\ker\bar\partial$. On the other hand, considering the $3$-space decomposition (\ref{eqn:BC-3sp-decomp}) of $C^{\infty}_{2, \, 0}(X, \C)$ induced by the Bott-Chern Laplacian $\Delta_{BC}:C^{\infty}_{2,\,0}(X,\,\C)\to C^{\infty}_{2,\,0}(X,\,\C)$ associated with the metric $\omega$, we see that the form $\rho^{2,\,0}$ with minimal $L^2_\omega$ norm satisfying equations $(i)$ and $(ii)$ of (\ref{eqn:H-S_condition_bis}) is the unique such form lying in the orthogonal complement of $\ker\partial\cap\ker\bar\partial = \ker\Delta_{BC} \oplus \mbox{Im}\,\partial\bar\partial$ in $C^{\infty}_{2, \, 0}(X, \C)$, which is $\mbox{Im}\,\partial^{\star}_\omega + \mbox{Im}\,\bar\partial^{\star}_\omega $.  

 For the proof of formula (\ref{eqn:Neumann_torsion_H-S}), see Lemma \ref{Lem:Neumann_del-0} with $v=-\partial\omega$.   \hfill  $\Box$

\begin{Obs}\label{Obs:torsion-formula_dim3} When $\mbox{dim}_\C X=3$, formula (\ref{eqn:Neumann_torsion_H-S}) for the $(2,\,0)$-torsion form $\rho_\omega^{2,\,0}$ of any Hermitian-symplectic metric $\omega$ simplifies to \begin{equation}\label{eqn:Neumann_torsion_H-S_dim3}\rho_\omega^{2,\,0} = -\Delta^{''-1}\bar\partial^\star(\partial\omega),\end{equation}

\noindent where $\Delta^{''-1} = \Delta^{''-1}_\omega$ is the Green operator of the $\bar\partial$-Laplacian $\Delta'' = \Delta''_\omega:=\bar\partial\bar\partial^\star + \bar\partial^\star\bar\partial$ induced by $\omega$ via $\bar\partial^\star = \bar\partial^\star_\omega$.  

\end{Obs}

\noindent {\it Proof.} It is a standard and easily-verified fact that on any compact complex $n$-dimensional manifold, any $\bar\partial$-closed $(n-1,\,0)$-form is $\partial$-closed. Now, the $(2,\,0)$-form $\rho^{2,\,0}$ satisfying $\partial\rho^{2,\,0}=0$ and $\bar\partial\rho^{2,\,0}=-\partial\omega$ (cf.\!\! (\ref{eqn:H-S_condition_bis})) is unique up to the addition of an arbitrary $(2,\,0)$-form $\zeta\in\ker\partial\cap\ker\bar\partial$. When $n=3$, $n-1=2$, so $\ker\partial\cap\ker\bar\partial = \ker\bar\partial$ in bidegree $(2,\,0)$. Therefore, $\rho_\omega^{2,\,0}\in\mbox{Im}\,\bar\partial^{\star}_\omega$, i.e. $\rho_\omega^{2,\,0} = \bar\partial^{\star}\xi$ for some $(2,\,1)$-form $\xi$. We get $\Delta''\rho_\omega^{2,\,0} = \bar\partial^{\star}\bar\partial( \bar\partial^{\star}\xi) = -\bar\partial^{\star}(\partial\omega)$. This is equivalent to (\ref{eqn:Neumann_torsion_H-S_dim3}). \hfill $\Box$

\vspace{3ex}

 If $\omega_0$ is a Hermitian-symplectic metric on $X$, any $C^{\infty}$ positive definite $(1,\,1)$-form $\omega$ lying in the Aeppli cohomology class of $\omega_0$ is a Hermitian-symplectic metric. Indeed, by $(i)$ and $(ii)$ of (\ref{eqn:H-S_condition_bis}), $\partial\omega_0 = -\bar\partial\rho_0^{2,\,0}$ for some $\partial$-closed $(2,\,0)$-form $\rho_0^{2,\,0}$ on $X$. Meanwhile, $\omega = \omega_0 + \partial\bar{u} + \bar\partial u$ for some $(1,\,0)$-form $u$, so $\partial\omega = \partial\omega_0 + \partial\bar\partial u = -\bar\partial(\rho_0^{2,\,0} + \partial u)$. Moreover, $\rho_0^{2,\,0} + \partial u$ is $\partial$-closed since $\rho_0^{2,\,0}$ is. Therefore, $\omega$ is Hermitian-symplectic (cf. $(i)$ and $(ii)$ of (\ref{eqn:H-S_condition_bis}) which characterise the H-S property).

 By a {\it Hermitian-symplectic (H-S) Aeppli class} $\{\omega\}_A\in H^{1,\,1}_A(X,\,\R)$ we shall mean a real Aeppli cohomology class of bidegree $(1,\,1)$ that contains an H-S metric $\omega$. We denote by $${\cal HS}_X:=\bigg\{\{\omega\}_A\in H_A^{1,\,1}(X,\,\R)\,\mid\,\omega \hspace{1ex}\mbox{is an H-S metric on}\hspace{1ex} X\bigg\}\subset H_A^{1,\,1}(X,\,\R)$$ the set of all such classes. Moreover, for every such class $\{\omega\}_A\in{\cal HS}_X$, we denote by $${\cal S}_{\{\omega\}}:=\bigg\{\omega + \partial\bar{u} + \bar\partial u\,\mid\,u\in C^\infty_{1,\,0}(X,\,\C) \hspace{1ex}\mbox{such that}\hspace{1ex} \omega + \partial\bar{u} + \bar\partial u>0\bigg\}\subset\{\omega\}_A\cap C^\infty_{1,\,1}(X,\,\R)$$ the set of all (necessarily H-S) metrics in $\{\omega\}_A$. The set ${\cal S}_{\{\omega\}}$ is an {\it open convex subset} of the real affine space $\{\omega\}_A\cap C^\infty_{1,\,1}(X,\,\R) = \{\omega + \partial\bar{u} + \bar\partial u\,\mid\,u\in C^\infty_{1,\,0}(X,\,\C)\}$.

\begin{Def}\label{Def:F_energy-functional_H-S} Let $X$ be a compact complex Hermitian-symplectic manifold with $\mbox{dim}_{\C}X=n$. For the Aeppli cohomology class $\{\omega_0\}_A\in{\cal HS}_X$ of any Hermitian-symplectic metric $\omega_0$, we define the following {\bf energy functional}:

\begin{equation}\label{eqn:F_energy-functional_H-S}F : {\cal S}_{\{\omega_0\}} \to [0,\,+\infty), \hspace{3ex} F(\omega) = \int\limits_X|\rho_\omega^{2,\,0}|^2_\omega\,dV_\omega = ||\rho_\omega^{2,\,0}||^2_\omega,  \end{equation}

\noindent where $\rho_{\omega}^{2,\,0}$ is the $(2,\,0)$-torsion form of the Hermitian-symplectic metric $\omega\in{\cal S}_{\{\omega_0\}}$ defined in Lemma and Definition \ref{Lem-Def:minimal_rho}, while $|\,\,\,|_\omega$ is the pointwise norm and $||\,\,\,||_\omega$ is the $L^2$ norm induced by $\omega$.

\end{Def}

The first trivial observation that justifies the introduction of the functional $F$ is the following.

\begin{Lem}\label{Lem:vanishing-F-Kaehler} Let $\{\omega_0\}_A\in{\cal HS}_X$ and $\omega\in{\cal S}_{\{\omega_0\}}$. Then, the following equivalence holds:

\begin{equation}\label{eqn:vanishing-F-Kaehler}\omega \hspace{1ex} \mbox{is K\"ahler} \iff F(\omega)=0.\end{equation}

\end{Lem} 

\noindent {\it Proof.} If $\omega$ is K\"ahler, $\partial\omega=0$ and the minimal $L^2$-norm solution of the equation $\bar\partial\rho=0$ vanishes. Thus $\rho_\omega^{2,\,0}=0$, hence $F(\omega)=0$. Conversely, if $F(\omega)=0$, then $\rho_\omega^{2,\,0}$ vanishes identically on $X$, hence $\partial\omega = -\bar\partial\rho_\omega^{2,\,0} = 0$, so $\omega$ is K\"ahler.  \hfill $\Box$

\vspace{2ex}

 We now compute the critical points of the energy functional $F$.

 Note that definition (\ref{eqn:F_energy-functional}) of $F$ translates to

\begin{equation}\label{eqn:F_energy-functional_bis} F(\omega) = \int\limits_X\rho_\omega^{2,\,0}\wedge\star\overline{\rho^{2,\,0}_\omega} = \int\limits_X\rho_\omega^{2,\,0}\wedge\rho_\omega^{0,\,2}\wedge\frac{\omega^{n-2}}{(n-2)!}.\end{equation}

\noindent Indeed, $\overline{\rho^{2,\,0}_\omega} = \rho_\omega^{0,\,2}$ is primitive since it is of bidegree $(0,\,2)$, so $\star\overline{\rho^{2,\,0}_\omega} =\overline{\rho^{2,\,0}_\omega}\wedge\omega^{n-2}/(n-2)!$ by (\ref{eqn:prim-form-star-formula-gen}).

 We now fix a Hermitian-symplectic metric $\omega$ on $X$ and we vary it in its Aeppli class along the path $\omega + t\gamma$, where $\gamma=\partial\bar{u} + \bar\partial u\in C^\infty_{1,\,1}(X,\,\R)$ is a fixed real $(1,\,1)$-form chosen to be Aeppli cohomologous to zero. Recall that the $(2,\,0)$-torsion form $\rho_{\omega}^{2,\,0}$ satisfies the condition $\bar\partial\rho_{\omega}^{2,\,0} = -\partial\omega$ and has minimal $L^2_\omega$-norm with this property. We get

\begin{equation}\label{eqn:rho_plus_t_u}\bar\partial(\rho_{\omega}^{2,\,0} + t\partial u) = -\partial(\omega + t\gamma),\end{equation}

\noindent although $\rho_{\omega}^{2,\,0} + t\partial u$ need not be of minimal $L^2_{\omega + t\gamma}$-norm with this property. For every $t\in\R$ close to $0$, we define the new functional: \begin{eqnarray}\label{eqn:F-tilde_def}\nonumber\widetilde{F}_t(\omega) & := & \int\limits_X|\rho_{\omega}^{2,\,0} + t\partial u|^2_{\omega + t\gamma}\,\frac{(\omega + t\gamma)^n}{n!} = \int\limits_X(\rho_{\omega}^{2,\,0} + t\partial u)\wedge\star_{\omega + t\gamma}\,(\overline{\rho_{\omega}^{2,\,0}} + t\,\bar\partial\bar{u}) \\
 & = & \int\limits_X(\rho_{\omega}^{2,\,0} + t\partial u)\wedge(\overline{\rho_{\omega}^{2,\,0}} + t\bar\partial\bar{u})\wedge\frac{(\omega + t\gamma)^{n-2}}{(n-2)!}.\end{eqnarray}

 The properties of $\widetilde{F}_t$ are summed up in the following statement.

\begin{Prop}\label{Prop:F-tilde_properties} $(i)$\, The two energy functionals are related by the inequality: \begin{equation}\label{eqn:F_energy_inequality}\widetilde{F}_t(\omega) \geq F(\omega + t\gamma) \hspace{2ex} \mbox{for all} \hspace{1ex} t\in\R \hspace{1ex} \mbox{close to} \hspace{1ex} 0.\end{equation}

\vspace{1ex}

 $(ii)$\, The differential at $\omega$ of $F$ is given by the formula: \begin{eqnarray}\label{eqn:differential_F}\nonumber(d_{\omega}F)(\gamma) = \frac{d}{dt}_{|t=0}\widetilde{F}_t(\omega) & = & -2\,\mbox{Re}\,\langle\langle u,\,\bar\partial^{\star}\omega\rangle\rangle_\omega + 2\,\mbox{Re}\,\int\limits_X u\wedge\rho_\omega^{2,\,0}\wedge\overline{\rho_\omega^{2,\,0}}\wedge\bar\partial\bigg(\frac{\omega^{n-3}}{(n-3)!}\bigg) \\
 & = & - \langle\langle\gamma\,,\omega\rangle\rangle + 2\,\mbox{Re}\,\int\limits_X u\wedge\rho_\omega^{2,\,0}\wedge\overline{\rho_\omega^{2,\,0}}\wedge\bar\partial\bigg(\frac{\omega^{n-3}}{(n-3)!}\bigg),\end{eqnarray}

\noindent for every $(1,\,1)$-form $\gamma = \partial\bar{u} + \bar\partial u$.

\end{Prop}

\noindent {\it Proof.} $(i)$\, If $t$ is sufficiently close to $0$, $\omega + t\gamma>0$, hence $\omega + t\gamma$ is a Hermitian-symplectic metric. By $(ii)$ in Lemma and Definition \ref{Lem-Def:minimal_rho}, we have $\bar\partial\rho_{\omega + t\gamma}^{2,\,0} = -\partial(\omega + t\gamma)$ and $\rho_{\omega + t\gamma}^{2,\,0}$ has minimal $L^2_{\omega+t\gamma}$-norm with this property. Since $\rho_{\omega}^{2,\,0} + t\partial u$ solves the same equation as $\rho_{\omega + t\gamma}^{2,\,0}$ (cf. (\ref{eqn:rho_plus_t_u})), we conclude that $$\widetilde{F}_t(\omega) \geq \int\limits_X|\rho_{\omega + t\gamma}^{2,\,0}|^2_{\omega + t\gamma}\,\frac{(\omega + t\gamma)^n}{n!} = F(\omega + t\gamma), \hspace{2ex} t\in\R \hspace{1ex} \mbox{close to} \hspace{1ex} 0.$$

$(ii)$\, Since $\widetilde{F}_t(\omega) - F(\omega + t\gamma)\geq 0$ for all $t\in\R$ close to $0$ and since $\widetilde{F}_0(\omega) = F(\omega)$, the smooth function $t\mapsto \widetilde{F}_t(\omega) - F(\omega + t\gamma)$ achieves a minimum at $t=0$. Hence, its derivative vanishes at $t=0$. We get: $$\frac{d}{dt}_{|t=0}\widetilde{F}_t(\omega) = \frac{d}{dt}_{|t=0} F(\omega + t\gamma) = (d_{\omega}F)(\gamma),$$

\noindent which is precisely the first identity in (\ref{eqn:differential_F}).

 We now prove the second identity in (\ref{eqn:differential_F}) starting from (\ref{eqn:F-tilde_def}). For all $t\in\R$ close to $0$, we get: \begin{eqnarray}\nonumber\frac{d}{dt}_{|t=0}\widetilde{F}_t(\omega) & = &  \frac{d}{dt}_{|t=0}\int\limits_X(\rho_{\omega}^{2,\,0} + t\partial u)\wedge(\overline{\rho_{\omega}^{2,\,0}} + t\bar\partial\bar{u})\wedge\frac{(\omega + t\gamma)^{n-2}}{(n-2)!} \\
\nonumber & = & \int\limits_X\partial u\wedge\overline{\rho_{\omega}^{2,\,0}}\wedge\frac{\omega^{n-2}}{(n-2)!} + \int\limits_X\bar\partial\bar{u}\wedge\rho_{\omega}^{2,\,0}\wedge\frac{\omega^{n-2}}{(n-2)!} \\
\nonumber & & + \int\limits_X\rho_{\omega}^{2,\,0}\wedge\overline{\rho_{\omega}^{2,\,0}}\wedge\frac{\omega^{n-3}}{(n-3)!}\wedge(\partial\bar{u} + \bar\partial u).\end{eqnarray}

\noindent Applying Stokes's theorem in each integral to remove the derivatives from $u$ and $\bar{u}$, we get: \begin{eqnarray}\nonumber\frac{d}{dt}_{|t=0}\widetilde{F}_t(\omega) & = &  \int\limits_X u\wedge\partial\overline{\rho_{\omega}^{2,\,0}}\wedge\frac{\omega^{n-2}}{(n-2)!} +\int\limits_X u\wedge\overline{\rho_{\omega}^{2,\,0}}\wedge\partial\bigg(\frac{\omega^{n-2}}{(n-2)!}\bigg) \\
\nonumber & & + \int\limits_X\bar{u}\wedge\bar\partial\rho_{\omega}^{2,\,0}\wedge\frac{\omega^{n-2}}{(n-2)!} + \int\limits_X\bar{u}\wedge\rho_{\omega}^{2,\,0}\wedge\bar\partial\bigg(\frac{\omega^{n-2}}{(n-2)!}\bigg)\\    \nonumber & & +\int\limits_X\bar{u}
\wedge\partial\rho_{\omega}^{2,\,0}\wedge\overline{\rho_{\omega}^{2,\,0}}\wedge\frac{\omega^{n-3}}{(n-3)!} +\int\limits_X\bar{u}\wedge\rho_{\omega}^{2,\,0}\wedge\partial\overline{\rho_{\omega}^{2,\,0}}\wedge\frac{\omega^{n-3}}{(n-3)!} \\
\nonumber & & + \int\limits_X\bar{u}\wedge\rho_{\omega}^{2,\,0}\wedge\overline{\rho_{\omega}^{2,\,0}}\wedge\partial\bigg(\frac{\omega^{n-3}}{(n-3)!}\bigg) + \int\limits_X u\wedge\rho_{\omega}^{2,\,0}\wedge\overline{\rho_{\omega}^{2,\,0}}\wedge\bar\partial\bigg(\frac{\omega^{n-3}}{(n-3)!}\bigg) \\
\nonumber &  & +\int\limits_X u\wedge\bar\partial\rho_{\omega}^{2,\,0}\wedge\overline{\rho_{\omega}^{2,\,0}}\wedge\frac{\omega^{n-3}}{(n-3)!} +\int\limits_X u\wedge\rho_{\omega}^{2,\,0}\wedge\bar\partial\overline{\rho_{\omega}^{2,\,0}}\wedge\frac{\omega^{n-3}}{(n-3)!}.\end{eqnarray}

\noindent Grouping the terms on the r.h.s. according to whether the integrands are divisible by $u$ or by $\bar{u}$ and using the identities $\partial\overline{\rho_\omega^{2,\,0}} = -\bar\partial\omega$ and $\bar\partial\rho_\omega^{2,\,0} = -\partial\omega$ , we get: \begin{eqnarray}\nonumber\frac{d}{dt}_{|t=0}\widetilde{F}_t(\omega) & = &  -\int\limits_X u\wedge\bigg[\bar\partial\omega\wedge\frac{\omega^{n-2}}{(n-2)!} + \bigg(-\overline{\rho_{\omega}^{2,\,0}}\wedge\partial\frac{\omega^{n-2}}{(n-2)!} + \partial\omega\wedge\overline{\rho_{\omega}^{2,\,0}}\wedge\frac{\omega^{n-3}}{(n-3)!}\bigg)\bigg] \\
\nonumber &  & + \int\limits_X u\wedge\bigg[\rho_{\omega}^{2,\,0}\wedge\bar\partial\overline{\rho_{\omega}^{2,\,0}}\wedge\frac{\omega^{n-3}}{(n-3)!} + \rho_{\omega}^{2,\,0}\wedge\overline{\rho_{\omega}^{2,\,0}}\wedge\bar\partial\frac{\omega^{n-3}}{(n-3)!}\bigg]  \\
\nonumber &  & -\int\limits_X \bar{u}\wedge\bigg[\partial\omega\wedge\frac{\omega^{n-2}}{(n-2)!} + \bigg(-\rho_{\omega}^{2,\,0}\wedge\bar\partial\frac{\omega^{n-2}}{(n-2)!} + \bar\partial\omega\wedge\rho_{\omega}^{2,\,0}\wedge\frac{\omega^{n-3}}{(n-3)!}\bigg)\bigg] \\ 
\nonumber &  & + \int\limits_X \bar{u}\wedge\bigg[\overline{\rho_{\omega}^{2,\,0}}\wedge\partial\rho_{\omega}^{2,\,0}\wedge\frac{\omega^{n-3}}{(n-3)!} + \rho_{\omega}^{2,\,0}\wedge\overline{\rho_{\omega}^{2,\,0}}\wedge\partial\frac{\omega^{n-3}}{(n-3)!}\bigg].\end{eqnarray}

\noindent Now, the terms on the first two lines on the r.h.s. above are respectively conjugated to the terms on the third and fourth lines, while the two inner large paratheses on lines $1$ and $3$ vanish since $\partial\omega^{n-2}/(n-2)! = \partial\omega\wedge\omega^{n-3}/(n-3)!$. On the other hand, we recall that $\partial\rho_\omega^{2,\,0}=0$, hence also $\bar\partial\overline{\rho_\omega^{2,\,0}}=0$. Thus, the two integrals containing these factors on the r.h.s. above vanish. We are reduced to \begin{eqnarray}\nonumber\frac{d}{dt}_{|t=0}\widetilde{F}_t(\omega) & = &  -\int\limits_X u\wedge\bigg[\bar\partial\frac{\omega^{n-1}}{(n-1)!} - \rho_{\omega}^{2,\,0}\wedge\overline{\rho_{\omega}^{2,\,0}}\wedge\bar\partial\frac{\omega^{n-3}}{(n-3)!}\bigg] \\
 &  & -\int\limits_X \bar{u}\wedge\bigg[\partial\frac{\omega^{n-1}}{(n-1)!} - \rho_{\omega}^{2,\,0}\wedge\overline{\rho_{\omega}^{2,\,0}}\wedge\partial\frac{\omega^{n-3}}{(n-3)!}\bigg],\end{eqnarray}

\noindent or equivalently, to \begin{eqnarray}\label{eqn:d_dt_F_tilde}\frac{d}{dt}_{|t=0}\widetilde{F}_t(\omega) = -2\,\mbox{Re}\,\int\limits_X u\wedge\bar\partial\bigg(\frac{\omega^{n-1}}{(n-1)!}\bigg) + 2\,\mbox{Re}\,\int\limits_X u\wedge\rho_\omega^{2,\,0}\wedge\overline{\rho_\omega^{2,\,0}}\wedge\bar\partial\bigg(\frac{\omega^{n-3}}{(n-3)!}\bigg).\end{eqnarray}

 Now, from $\star_\omega\omega = \omega^{n-1}/(n-1)!$ and $\partial^{\star} = -\star\bar\partial\star$, we get: $$\bar\partial\bigg(\frac{\omega^{n-1}}{(n-1)!}\bigg) = \bar\partial\star\omega = \star(-\star\bar\partial\star)\,\omega = \star\partial^{\star}\omega = \star\,\overline{\bar\partial^{\star}\omega},$$

\noindent hence \begin{equation}\label{eqn:u_wedge_inner-product}u\wedge \bar\partial\bigg(\frac{\omega^{n-1}}{(n-1)!}\bigg) = \langle u,\,\bar\partial^{\star}\omega\rangle_\omega \, dV_\omega.\end{equation}

 Thus, (\ref{eqn:d_dt_F_tilde}) and (\ref{eqn:u_wedge_inner-product}) prove the second identity in (\ref{eqn:differential_F}). The third identity in (\ref{eqn:differential_F}) is obvious.  \hfill $\Box$

\begin{Cor}\label{Cor:critical-points_energy_n=3} Suppose $n=3$. Then a Hermitian-symplectic metric $\omega$ on a compact complex manifold $X$ of dimension $3$ is a {\bf critical point} of the energy functional $F$ {\bf if and only if} $\omega$ is {\bf K\"ahler}.

\end{Cor}

\noindent {\it Proof.} It is obvious that every K\"ahler metric $\omega$ is a critical point for $F$ since $\partial\omega=0$, hence $\rho_\omega^{2,\,0} = 0$. 

 If $n=3$, $\bar\partial\omega^{n-3}=0$, so (\ref{eqn:differential_F}) reduces to $(d_\omega F)(\gamma) = -2\,\mbox{Re}\,\langle\langle u,\,\bar\partial^{\star}\omega\rangle\rangle_\omega$.

 Now, a metric $\omega$ is a critical point of $F$ if and only if $(d_\omega F)(\gamma) = 0$ for every $\gamma = \partial\bar{u} + \bar\partial u$. By the above discussion, this amounts to $\mbox{Re}\,\langle\langle u,\,\bar\partial^{\star}\omega\rangle\rangle_\omega = 0$ for every $(1,\,0)$-form $u$. Thus, if $\omega$ is a critical point of $F$, by taking $u=\bar\partial^{\star}\omega$ we get $\bar\partial^{\star}\omega = 0$. This is equivalent to $\omega$ being balanced. However, $\omega$ is already SKT since it is Hermitian-symplectic, so $\omega$ must be K\"ahler by Proposition \ref{Prop:SKT+bal}.    \hfill $\Box$

\vspace{3ex}

\begin{Cor}\label{Cor:energy_M-A_mass} Let $X$ be a compact complex manifold of dimension $n=3$ admitting Hermitian-symplectic metrics. Then, for every Aeppli-cohomologous Hermitian-symplectic metrics $\omega$ and $\omega_\eta$: \begin{equation}\label{eqn:A_cohomologous_H-S}\omega_\eta = \omega + \partial\bar\eta + \bar\partial\eta >0,  \hspace{3ex} \mbox{with} \hspace{1ex} \eta\in C^{\infty}_{1,\,0}(X,\,\C),\end{equation}

\noindent the respective $(2,\,0)$-torsion forms $\rho_\omega^{2,\,0}$ and $\rho_\eta^{2,\,0}:=\rho_{\omega_\eta}^{2,\,0}$ satisfy the identity: \begin{equation}\label{eqn:torsion_A_cohomologous_H-S1}||\rho_\eta^{2,\,0}||^2_{\omega_\eta} + \int\limits_X\frac{\omega_\eta^3}{3!} = ||\rho_\omega^{2,\,0}||^2_\omega + \int\limits_X\frac{\omega^3}{3!}\end{equation}

\noindent and are related by \begin{equation}\label{eqn:torsion_A_cohomologous_H-S2}\rho_\eta^{2,\,0} = \rho_\omega^{2,\,0} + \partial\eta.\end{equation}

 In particular, if $\partial\eta=0$ (a condition that is equivalent to $\omega_\eta-\omega$ being $d$-exact), we are reduced to $\rho_\eta^{2,\,0} = \rho_\omega^{2,\,0}$ and \begin{equation}\label{eqn:torsion_A_cohomologous_H-S_bis}||\rho_\omega^{2,\,0}||^2_{\omega_\eta} = ||\rho_\omega^{2,\,0}||^2_\omega + \int\limits_X \rho_\omega^{2,\,0}\wedge\overline{\rho_\omega^{2,\,0}}\wedge(\omega_\eta - \omega).\end{equation}

\end{Cor}

\noindent{\it Proof.} In arbitrary dimension $n$, we compute the differential of the map $${\cal S}_{\{\omega_0\}}\ni\omega\mapsto\int\limits_X\frac{\omega^n}{n!}:=\mbox{Vol}_\omega(X)$$

\noindent when the metric $\omega$ varies in its Aeppli cohomology class $\{\omega_0\}_A$. For any real, Aeppli null-cohomologous $(1,\,1)$-form $\gamma = \partial\bar{u} + \bar\partial u$ (with $u\in C^{\infty}_{1,\,0}(X,\,\C)$), we have \begin{eqnarray}\nonumber\frac{d}{dt}_{|t=0}\int\limits_X\frac{(\omega + t\gamma)^n}{n!} & = & \frac{1}{(n-1)!}\,\int\limits_X\omega^{n-1}\wedge\gamma = 2\,\mbox{Re}\,\int\limits_X\bar\partial u \wedge\frac{\omega^{n-1}}{(n-1)!} = 2\,\mbox{Re}\,\int\limits_X u \wedge\bar\partial\star\omega \\
\nonumber & = & 2\,\mbox{Re}\,\int\limits_X u \wedge\star\bigg(-\star\bar\partial\star\omega\bigg) = 2\,\mbox{Re}\,\int\limits_X u \wedge\star\partial^{\star}\omega = 2\,\mbox{Re}\,\int\limits_X u \wedge\star\overline{\bar\partial^{\star}\omega} \\
\nonumber & = &  2\,\mbox{Re}\,\langle\langle u,\, \bar\partial^{\star}\omega\rangle\rangle.\end{eqnarray}

\noindent Together with (\ref{eqn:differential_F}) (recall that $n=3$ here), this identity shows that the differential at $\omega$ of the map $${\cal S}_{\{\omega_0\}}\ni\omega\mapsto ||\rho_\omega^{2,\,0}||^2_\omega + \int\limits_X\frac{\omega^3}{3!}$$

\noindent vanishes identically. Therefore, this map is constant on the Hermitian-symplectic metrics lying in a same Aeppli cohomology class $\{\omega_0\}_A$. This proves (\ref{eqn:torsion_A_cohomologous_H-S1}). 

 To prove (\ref{eqn:torsion_A_cohomologous_H-S2}), recall that definition (\ref{eqn:H-S_condition_bis}) of the $(2,\,0)$-torsion forms implies the following relations: \begin{equation}\label{eqn:torsion-forms_comparison_def} (i)\,\bar\partial(\rho_\eta^{2,\,0} - \partial\eta) = -\partial\omega \hspace{2ex} \mbox{and} \hspace{2ex} (ii)\, ||\rho_\eta^{2,\,0} - \partial\eta||_\omega \geq ||\rho_\omega^{2,\,0}||_\omega,\end{equation}

\noindent where $(ii)$ follows from $(i)$ and from the $L^2_\omega$-norm minimality of $\rho_\omega^{2,\,0}$ among the $(2,\,0)$-forms $\rho$ solving the equation $\bar\partial\rho = -\partial\omega$.  

 Now, (\ref{eqn:torsion_A_cohomologous_H-S1}) gives the first of the following identities: \begin{eqnarray}\label{eqn:torsion_L2_M-A_computation1}||\rho_\omega^{2,\,0}||^2_\omega + \int\limits_X\frac{\omega^3}{3!} & = & ||\rho_\eta^{2,\,0}||^2_{\omega_\eta} + \int\limits_X\frac{\omega_\eta^3}{3!} = \int\limits_X \rho_\eta^{2,\,0}\wedge\overline{\rho_\eta^{2,\,0}}\wedge\omega_\eta + \int\limits_X\frac{\omega_\eta^3}{3!}.\end{eqnarray}

\noindent On the other hand, we have: \begin{eqnarray}\label{eqn:torsion_L2_M-A_computation2}\nonumber(\omega_\eta + \rho_\eta^{2,\,0} + \overline{\rho_\eta^{2,\,0}})^3 & = & \omega_\eta^3 + 3\,\omega_\eta^2\wedge(\rho_\eta^{2,\,0} + \overline{\rho_\eta^{2,\,0}}) + 3\,\omega_\eta\wedge(\rho_\eta^{2,\,0} + \overline{\rho_\eta^{2,\,0}})^2 + (\rho_\eta^{2,\,0} + \overline{\rho_\eta^{2,\,0}})^3  \\
 & = & \omega_\eta^3 + 6\,\omega_\eta\wedge\rho_\eta^{2,\,0}\wedge\overline{\rho_\eta^{2,\,0}},\end{eqnarray}    

\noindent where the last identity follows from the cancellation of several terms for bidegree reasons. Putting (\ref{eqn:torsion_L2_M-A_computation1}) and (\ref{eqn:torsion_L2_M-A_computation2}) together, we get: \begin{eqnarray}\label{eqn:torsion_L2_M-A_computation3}\nonumber ||\rho_\omega^{2,\,0}||^2_\omega + \frac{1}{3!}\,\int\limits_X\omega^3 & = & \frac{1}{3!}\,\int\limits_X(\omega_\eta + \rho_\eta^{2,\,0} + \overline{\rho_\eta^{2,\,0}})^3 =  \frac{1}{3!}\, \int\limits_X[\omega + (\rho_\eta^{2,\,0} - \partial\eta) + (\overline{\rho_\eta^{2,\,0}} - \bar\partial\bar\eta) + d(\eta + \bar\eta)]^3 \\
\nonumber & \stackrel{(a)}{=} & \frac{1}{3!}\, \int\limits_X[\omega + (\rho_\eta^{2,\,0} - \partial\eta) + (\overline{\rho_\eta^{2,\,0}} - \bar\partial\bar\eta)]^3 \stackrel{(b)}{=} \frac{1}{3!}\, \int\limits_X\omega^3 + \int\limits_X(\rho_\eta^{2,\,0} - \partial\eta)\wedge(\overline{\rho_\eta^{2,\,0}} - \bar\partial\eta)\wedge\omega  \\
& = & ||\rho_\eta^{2,\,0} - \partial\eta||^2_\omega + \frac{1}{3!}\, \int\limits_X\omega^3 \stackrel{(c)}{\geq} ||\rho_\omega^{2,\,0}||^2_\omega + \frac{1}{3!}\, \int\limits_X\omega^3.\end{eqnarray}

$\cdot$ Identity $(a)$ followed from Stokes's theorem and the $d$-closedness of the form $\omega + (\rho_\eta^{2,\,0} - \partial\eta) + (\overline{\rho_\eta^{2,\,0}} - \bar\partial\bar\eta)$ that is seen through the following very simple computation: \begin{eqnarray}\nonumber d[\omega + (\rho_\eta^{2,\,0} - \partial\eta) + (\overline{\rho_\eta^{2,\,0}} - \bar\partial\bar\eta)] &  = & \partial\omega + \bar\partial\omega + \bar\partial\rho_\eta^{2,\,0} - \bar\partial\partial\eta + \partial\overline{\rho_\eta^{2,\,0}} - \partial\bar\partial\bar\eta \\
\nonumber & = & \partial\omega + \bar\partial\omega -\partial(\omega + \partial\bar\eta + \bar\partial\eta) - \bar\partial\partial\eta - \bar\partial(\omega + \bar\partial\eta + \partial\bar\eta) - \partial\bar\partial\bar\eta \\
\nonumber & = & - (\partial\bar\partial\eta + \bar\partial\partial\eta) - (\bar\partial\partial\bar\eta + \partial\bar\partial\bar\eta) = 0,\end{eqnarray}

\noindent where the second identity followed from $\bar\partial\rho_\eta^{2,\,0} = -\partial\omega_\eta = -\partial(\omega + \bar\partial\eta + \partial\bar\eta)$ and from the conjugated expression.

$\cdot$ Identity $(b)$ in (\ref{eqn:torsion_L2_M-A_computation3}) followed from the analogue of (\ref{eqn:torsion_L2_M-A_computation2}) in this context, while inequality $(c)$ in (\ref{eqn:torsion_L2_M-A_computation3}) followed from part $(ii)$ of (\ref{eqn:torsion-forms_comparison_def}).

 We see that the first and the last terms in (\ref{eqn:torsion_L2_M-A_computation3}) are equal. This forces $(c)$ to be an equality, hence part $(ii)$ of (\ref{eqn:torsion-forms_comparison_def}) must be an equality. This means that $\rho_\eta^{2,\,0}-\partial\eta$ and $\rho_\omega^{2,\,0}$ are both solutions of the equation $\bar\partial\rho^{2,\,0} = -\partial\omega$ (see part $(i)$ of (\ref{eqn:torsion-forms_comparison_def}) and part $(ii)$ of (\ref{eqn:H-S_condition_bis})) and have {\it equal} $L^2_\omega$-norms. Since $\rho_\omega^{2,\,0}$ is the minimal $L^2_\omega$-norm solution, we infer that $\rho_\eta^{2,\,0}-\partial\eta = \rho_\omega^{2,\,0}$ by the uniqueness of the minimal $L^2_\omega$-norm solution. This proves (\ref{eqn:torsion_A_cohomologous_H-S2}). 

 Finally, we write $\omega_\eta = \omega + (\omega_\eta-\omega)$ and $$||\rho_\omega^{2,\,0}||^2_{\omega_\eta} = \int\limits_X\rho_\omega^{2,\,0}\wedge\overline{\rho_\omega^{2,\,0}}\wedge\omega_\eta = ||\rho_\omega^{2,\,0}||^2_\omega + \int\limits_X\rho_\omega^{2,\,0}\wedge\overline{\rho_\omega^{2,\,0}}\wedge(\omega_\eta-\omega).$$

\noindent This is (\ref{eqn:torsion_A_cohomologous_H-S_bis}). \hfill $\Box$

\vspace{3ex}

The main takeaway from Corollary \ref{Cor:energy_M-A_mass} is that the sum $F(\omega) + \mbox{Vol}_\omega(X)$ (where $\mbox{Vol}_\omega(X):=\int_X\omega^3/3!$) remains {\bf constant} when $\omega$ ranges over the (necessarily Hermitian-symplectic) metrics in the Aeppli cohomology class of a fixed Hermitian-symplectic metric $\omega_0$. This invariant attached to any Aeppli class of Hermitian-symplectic metrics generalises the classical volume of a K\"ahler class and constitutes one of our main findings in this work.

\begin{Def}\label{Def:A-invariant} Let $X$ be a $3$-dimensional compact complex manifold supposed to carry Hermitian-symplectic metrics. For any such metric $\omega$ on $X$, the constant \begin{equation}\label{eqn:A-invariant}A=A_{\{\omega\}_A}:= F(\omega) + \mbox{Vol}_\omega(X)>0\end{equation} depending only on $\{\omega\}_A$ is called the {\bf generalised volume} of the Hermitian-symplectic Aeppli class $\{\omega\}_A$.

\end{Def}

\subsection{Case of SKT metrics on $\partial\bar\partial$-manifolds}\label{subsection:SKT_ddbar}

In this subsection, we discuss an analogous functional in the special case of a $\partial\bar\partial$-manifold admitting SKT metrics.

Note that any Hermitian-symplectic metric is SKT on any manifold $X$. Moreover, for every Hermitian-symplectic metric $\omega$ on $X$, the set of all Hermitian-symplectic metrics in the Aeppli class $\{\omega\}_A$ coincides with the set of all SKT metrics in the Aeppli class $\{\omega\}_A$. Conversely, SKT manifolds need not be Hermitian-symplectic, but on $\partial\bar\partial$-manifolds, the notions of H-S and SKT metrics coincide.

\begin{Lem}\label{Lem:torsion_2-0_form_SKT} For every SKT metric $\omega$ on a compact $\partial\bar\partial$-manifold $X$, there exists a unique smooth $(2,\,0)$-form $\Gamma_\omega$ on $X$ such that \begin{equation}\label{eqn:torsion_2-0_form}(i)\,\,\bar\partial\Gamma_\omega = -\partial\omega  \hspace{2ex} \mbox{and} \hspace{2ex} (ii)\,\,\Gamma_\omega\in\mbox{Im}\,\bar\partial^{\star}_\omega,\end{equation}

\noindent where the subscript $\omega$ indicates that the formal adjoint is computed w.r.t. the $L^2$ inner product defined by $\omega$. 

 The form $\Gamma_\omega$ will be called the {\bf $(2,\,0)$-torsion form} of the SKT metric $\omega$. It is given by the von Neumann-type formula: \begin{equation}\label{eqn:torsion_2-0_form_formula}\Gamma_\omega = -\Delta^{''-1}\bar\partial^{\star}(\partial\omega),\end{equation}

\noindent where $\Delta^{''-1}$ is the Green operator of the Laplacian $\Delta'' = \bar\partial\bar\partial^{\star} + \bar\partial^{\star}\bar\partial$ induced by the metric $\omega$.

\end{Lem}

\noindent {\it Proof.} The $(2,\,1)$-form $\partial\omega$ is $d$-closed (thanks to the SKT assumption on $\omega$) and $\partial$-exact, hence by the $\partial\bar\partial$-assumption on $X$ it is also $\bar\partial$-exact. This means that the equation $\bar\partial\Gamma = -\partial\omega$ is solvable. Its solutions $\Gamma$ are unique up to the addition of any element in $\ker\bar\partial$, so the minimal $L^2_\omega$-norm solution is the unique solution lying in the orthogonal complement of $\ker\bar\partial$, which is $\mbox{Im}\,\bar\partial^{\star}_\omega$. The von Neumann formula is well known and can be easily proved: $\bar\partial(-\Delta^{''-1}\bar\partial^{\star}(\partial\omega)) = -\partial\omega$ (immediate verification) and $\Delta^{''-1}\bar\partial^{\star}(\partial\omega) = \bar\partial^{\star}\Delta^{''-1}(\partial\omega)\in\mbox{Im}\,\bar\partial^{\star}$.  \hfill $\Box$

\vspace{2ex}

 The following is a very simple observation.

\begin{Lem}\label{Lem:Gamma_del-closed} The $(2,\,0)$-torsion form $\Gamma_\omega$ of any SKT metric $\omega$ on a compact $\partial\bar\partial$-manifold $X$ has the property: \begin{equation}\label{eqn:Gamma_del-closed}\partial\Gamma_\omega = 0.\end{equation}

\end{Lem}

\noindent {\it Proof.} The $(3,\,0)$-form $\partial\Gamma_\omega$ is $\partial$-exact (obviously) and $d$-closed (since $\bar\partial(\partial\Gamma_\omega) = -\partial(\bar\partial\Gamma_\omega) = \partial^2\omega = 0$), hence it must be $\partial\bar\partial$-exact thanks to the $\partial\bar\partial$ assumption on $X$. This means that there exists a $(2,\,-1)$-form $\zeta$ (which must vanish for bidegree reasons) such that $\partial\bar\partial\zeta = \partial\Gamma_\omega$. Then $\partial\Gamma_\omega$ vanishes since $\zeta = 0$.  \hfill $\Box$

\vspace{2ex}

 We now define a new energy functional by the $L^2$-norm of the $(2,\,0)$-torsion form $\Gamma_\omega$.

\begin{Def}\label{Def:F_energy-functional} Let $X$ be a compact SKT $\partial\bar\partial$-manifold with $\mbox{dim}_{\C}X=n$. For every Aeppli cohomology class $\{\omega_0\}_A$ representable by an SKT metric, we define the following {\bf energy functional}: \begin{equation}\label{eqn:F_energy-functional}F : {\cal S}_{\{\omega_0\}} \to [0,\,+\infty), \hspace{3ex} F(\omega) = \int\limits_X|\Gamma_\omega|^2_\omega\,dV_\omega = ||\Gamma_\omega||^2_\omega,\end{equation}

\noindent where $\Gamma_{\omega}$ is the $(2,\,0)$-torsion form of the SKT metric $\omega\in{\cal S}_{\{\omega_0\}}$ defined in Lemma \ref{Lem:torsion_2-0_form_SKT}.

\end{Def}

The remaining arguments are identical to those given in $\S.$\ref{subsection:H-S_n=3_critical-points} if we replace $\rho_\omega^{2,\,0}$ with $\Gamma_\omega$. Recall that by (\ref{eqn:Gamma_del-closed}) we have $\partial\Gamma_\omega=0$ (cf. $(i)$ of (\ref{eqn:H-S_condition_bis})). 

 The first variation of $F$ can be computed as in $\S.$\ref{subsection:H-S_n=3_critical-points}. We get \begin{Prop}\label{Prop:F-tilde_properties_SKT} The differential of $F$ at any SKT metric $\omega$ is given by the formula: \begin{eqnarray}\label{eqn:differential_F_H-S}\nonumber(d_{\omega}F)(\gamma) & = & -2\,\mbox{Re}\,\langle\langle u,\,\bar\partial^{\star}\omega\rangle\rangle_\omega + 2\,\mbox{Re}\,\int\limits_X u\wedge\Gamma_\omega\wedge\overline{\Gamma}_\omega\wedge\bar\partial\bigg(\frac{\omega^{n-3}}{(n-3)!}\bigg) \\
 & = & - \langle\langle\gamma\,,\omega\rangle\rangle + 2\,\mbox{Re}\,\int\limits_X u\wedge\Gamma_\omega\wedge\overline{\Gamma}_\omega\wedge\bar\partial\bigg(\frac{\omega^{n-3}}{(n-3)!}\bigg)  \end{eqnarray}

\noindent for every $(1,\,1)$-form $\gamma = \partial\bar{u} + \bar\partial u$.

 In particular, {\bf if $n=3$}, an {SKT} metric $\omega$ on $X$ is a {\bf critical point} of the energy functional $F$ {\bf if and only if} $\omega$ is {\bf K\"ahler}.

\end{Prop}

\subsection{Variation of the $(2,\,0)$-torsion form for $\partial\bar\partial$-cohomologous metrics}\label{subsection:variation_torsion}

We first show that the $(2,\,0)$-torsion form of a Hermitian-symplectic metric does not change when the metric changes only by an element in $\mbox{Im}\,\partial\bar\partial$. The next statement can be compared with Corollary \ref{Cor:energy_M-A_mass}: it supposes more and achieves more.

\begin{Prop}\label{Prop:torsion_variation_BC} Let $X$ be a compact complex manifold with $\mbox{dim}_\C X=3$. Suppose that $\omega>0$ and $\widetilde\omega = \omega + i\partial\bar\partial\varphi>0$ are {\bf SKT} metrics on $X$.

 $(i)$\, For every form $\rho^{2,\,0}\in C^{\infty}_{2,\,0}(X,\,\C)$ such that $\partial\rho^{2,\,0}=0$ and $\bar\partial\rho^{2,\,0} = -\partial\omega$, the $L^2$-norms of $\rho^{2,\,0}$ w.r.t. $\widetilde\omega$ and $\omega$ are related in the following way: \begin{equation}\label{eqn:torsion_variation_BC_L2}||\rho^{2,\,0}||^2_{\widetilde\omega} = ||\rho^{2,\,0}||^2_\omega - \frac{1}{2}\,\int\limits_X(\widetilde\omega - \omega)\wedge\omega^2.\end{equation}

\noindent This relation is equivalent to \begin{equation}\label{eqn:torsion_variation_BC_L2_bis}||\rho^{2,\,0}||^2_{\widetilde\omega} + \int\limits_X\frac{\widetilde\omega^3}{3!} = ||\rho^{2,\,0}||^2_\omega + \int\limits_X\frac{\omega^3}{3!}.\end{equation}

 $(ii)$\, If $\omega>0$ and $\widetilde\omega = \omega + i\partial\bar\partial\varphi>0$ are {\bf Hermitian-symplectic} metrics, their $(2,\,0)$-torsion forms coincide, i.e. \begin{equation}\label{eqn:torsion_variation_BC}\rho^{2,\,0}_{\widetilde\omega} = \rho^{2,\,0}_\omega.\end{equation}

\end{Prop}

\noindent {\it Proof.} $(i)$\, From the assumptions, we get the following identities: \begin{eqnarray}\label{eqn:L2_norms_torsions_computations}\nonumber ||\rho^{2,\,0}||^2_{\widetilde\omega} & = & \int\limits_X\rho^{2,\,0}\wedge\overline{\rho^{2,\,0}}\wedge\widetilde\omega = \int\limits_X\rho^{2,\,0}\wedge\overline{\rho^{2,\,0}}\wedge\omega + \int\limits_X\rho^{2,\,0}\wedge\overline{\rho^{2,\,0}}\wedge i\partial\bar\partial\varphi \\
 \nonumber & \stackrel{(a)}{=} & ||\rho^{2,\,0}||^2_\omega -i\,\int\limits_X\rho^{2,\,0}\wedge\partial\overline{\rho^{2,\,0}}\wedge\bar\partial\varphi \\
 & \stackrel{(b)}{=} & ||\rho^{2,\,0}||^2_\omega -i\,\int\limits_X\varphi\, \bar\partial\rho^{2,\,0}\wedge\partial\overline{\rho^{2,\,0}}  \stackrel{(c)}{=} ||\rho^{2,\,0}||^2_\omega -i\,\int\limits_X\varphi\,\partial\omega\wedge\bar\partial\omega,\end{eqnarray} 

\noindent where $(a)$ and $(b)$ follow from Stokes combined with the identity $\partial\rho^{2,\,0}=0$ and its conjugate $\bar\partial\overline{\rho^{2,\,0}}=0$, while $(c)$ follows from the identity $\bar\partial\rho^{2,\,0} = -\partial\omega$ and its conjugate $\partial\overline{\rho^{2,\,0}} = -\bar\partial\omega$.

 Now, the SKT property of $\omega$ implies that $\partial\omega\wedge\bar\partial\omega = \partial(\omega\wedge\bar\partial\omega) = (1/2)\,\partial\bar\partial\omega^2$, so two further applications of Stokes yield the second identity below:

\begin{eqnarray}\label{eqn:L2_norms_torsions_computations_1}i\,\int\limits_X\varphi\,\partial\omega\wedge\bar\partial\omega = \frac{i}{2}\,\int\limits_X\varphi\,\partial\bar\partial\omega^2 = \frac{1}{2}\,\int\limits_X i\partial\bar\partial\varphi\wedge\omega^2 = \frac{1}{2}\,\int\limits_X (\widetilde\omega - \omega)\wedge\omega^2.\end{eqnarray}

 We now see that (\ref{eqn:L2_norms_torsions_computations}) and (\ref{eqn:L2_norms_torsions_computations_1}) prove (\ref{eqn:torsion_variation_BC_L2}) between them.  

 To prove the equivalence of (\ref{eqn:torsion_variation_BC_L2_bis}) and (\ref{eqn:torsion_variation_BC_L2}), we have to show that $$(1/6)\,\int_X(\widetilde\omega^3 - \omega^3) = (1/2)\,\int_X(\widetilde\omega - \omega)\wedge\omega^2.$$ Now, since $\widetilde\omega^2 = \omega^2 + 2\,i\partial\bar\partial\varphi\wedge\omega + (i\partial\bar\partial\varphi)^2$, we get: \begin{eqnarray}\nonumber \frac{1}{6}\,\int\limits_X(\widetilde\omega^3 - \omega^3) & = & \frac{1}{6}\,\int\limits_X(\widetilde\omega - \omega)\wedge(\widetilde\omega^2 + \widetilde\omega\wedge\omega + \omega^2) \\
\nonumber & = & \frac{1}{6}\,\int\limits_X(\widetilde\omega - \omega)\wedge\omega^2 + \frac{1}{3}\,\int\limits_X(\widetilde\omega - \omega)\wedge i\partial\bar\partial\varphi\wedge\omega + \frac{1}{6}\,\int\limits_X(\widetilde\omega - \omega)\wedge(i\partial\bar\partial\varphi)^2  \\
\nonumber & + & \frac{1}{6}\,\int\limits_X(\widetilde\omega - \omega)\wedge\omega^2 + \frac{1}{6}\,\int\limits_X(\widetilde\omega - \omega)\wedge i\partial\bar\partial\varphi \wedge\omega + \frac{1}{6}\,\int\limits_X(\widetilde\omega - \omega)\wedge\omega^2 \\
 \nonumber & = & 3\cdot \frac{1}{6}\,\int\limits_X(\widetilde\omega - \omega)\wedge\omega^2,\end{eqnarray}

\noindent since all the other terms vanish by Stokes, the identities $\partial(\widetilde\omega - \omega) = 0$ and $\bar\partial(\widetilde\omega - \omega) = 0$ and the SKT assumption on $\omega$.

$(ii)$\, The stronger H-S assumption on $\widetilde\omega$ and $\omega$ is only made to ensure the existence of the $(2,\,0)$-torsion forms $\rho^{2,\,0}_{\widetilde\omega}$ and $\rho^{2,\,0}_\omega$. The assumption $\widetilde\omega = \omega +i\partial\bar\partial\varphi$ implies $\partial\widetilde\omega = \partial\omega$, so $\rho^{2,\,0}_{\widetilde\omega}$ and $\rho^{2,\,0}_\omega$ are the minimal $L^2_{\widetilde\omega}$-norm solution, resp. the minimal $L^2_\omega$-norm solution, of the same equation $\bar\partial\rho^{2,\,0} = -\partial\omega$.

However, (\ref{eqn:torsion_variation_BC_L2}) shows that when $\rho^{2,\,0}$ ranges over the set of smooth $(2,\,0)$-forms $\rho^{2,\,0}$ satisfying the conditions $\partial\rho^{2,\,0}=0$ and $\bar\partial\rho^{2,\,0} = -\partial\omega$ (both of which are satisfied by both $\rho^{2,\,0}_{\widetilde\omega}$ and $\rho^{2,\,0}_\omega$), $||\rho^{2,\,0}||_{\widetilde\omega}$ is minimal if and only if $||\rho^{2,\,0}||_\omega$ is minimal since the discrepancy term 

$$-\frac{1}{2}\,\int\limits_X(\widetilde\omega - \omega)\wedge\omega^2$$

\noindent is independent of $\rho^{2,\,0}$ (depending only on the given metrics $\widetilde\omega$ and $\omega$). This means that the same $\rho^{2,\,0}$ achieves the minimal $L^2$ norm w.r.t. either of the metrics $\widetilde\omega$ and $\omega$. By uniqueness of the minimal $L^2$-norm solution of the $\bar\partial$ equation, we get $\rho^{2,\,0}_{\widetilde\omega} = \rho^{2,\,0}_\omega$.  \hfill $\Box$

\subsection{Volume form and Monge-Amp\`ere-type equation associated with an H-S metric}\label{subsection:vol-M-A_eq_H-S}

 We now digress briefly to point out another possible future use of the new invariant defined by the generalised volume. In fact, a new volume form that seems better suited to featuring in the right-hand term of complex Monge-Amp\`ere equations can be associated with every Hermitian-symplectic metric on a $3$-dimensional compact complex manifold.

\begin{Def}\label{Def:Introd_tilde-volume-form} If $\omega$ is a Hermitian-symplectic metric on a compact complex manifold $X$ with $\mbox{dim}_\C X=3$ and $\rho_\omega^{2,\,0}$ is the  $(2,\,0)$-torsion form of $\omega$, we define the following volume form on $X$: $$d\widetilde{V}_\omega := (1 + |\rho_\omega^{2,\,0}|^2_\omega)\,dV_\omega.$$

\end{Def}  

  The main interest in this volume form stems from the fact that its volume is independent of the choice of metric in a given Hermitian-symplectic Aeppli class, as follows from Corollary \ref{Cor:energy_M-A_mass}: \begin{eqnarray*}\int\limits_X d\widetilde{V}_{\omega_1} = \int\limits_X d\widetilde{V}_{\omega_2} = A, \hspace{3ex} \mbox{for all metrics}\hspace{1ex} \omega_1,\omega_2\in\{\omega\}_A,\end{eqnarray*} where $A = A_{\{\omega\}_A}>0$ is the {\it generalised volume} of the H-S Aeppli class $\{\omega\}_A$.

  Now, if $\omega$ is a Hermitian-symplectic metric on a manifold $X$ as above, it seems natural to consider the  Monge-Amp\`ere equation $$\frac{(\omega + i\partial\bar\partial\varphi)^3}{3!} = b\,d\widetilde{V}_{\omega},$$ subject to the condition $\omega + i\partial\bar\partial\varphi>0$, where $b>0$ is a given constant. By [TW10, Corollary 1], there exists a unique $b$ such that this equation is solvable. Moreover, for that $b$, the solution $\omega + i\partial\bar\partial\varphi>0$ is unique. Note that $$b= \frac{\mbox{Vol}_{\omega + i\partial\bar\partial\varphi}(X)}{A_{\{\omega\}_A}}\in(0,\,1]$$ since $A_{\{\omega\}_A} = F(\omega + i\partial\bar\partial\varphi) + \mbox{Vol}_{\omega + i\partial\bar\partial\varphi}(X)\geq\mbox{Vol}_{\omega + i\partial\bar\partial\varphi}(X)$. We hope that this can shed some light on the mysterious constant $b$ in this context.

  \section{Obstruction and estimates}\label{section:obs-est}

  In this section, we point out an obstruction to the Aeppli cohomology class of a given Hermitian-symplectic metric containing a K\"ahler metric.

 We start by observing that a class in the vector space $E_2^{0,\,2} = E_2^{0,\,2}(X)$ featuring on the second page of the Fr\"olicher spectral sequence of $X$ can be uniquely associated with every Hermitian-symplectic metric on $X$ and, in dimension $3$, even with every Aeppli cohomology class of such metrics. (Cf. [PSUb, Proposition 6.2.].)

 \begin{Lem-Def}\label{Def:E_2_H-S} Suppose that $\omega$ is a {\bf Hermitian-symplectic metric} on a compact complex $n$-dimensional manifold $X$.

   (i)\, The $(0,\,2)$-torsion form $\rho_\omega^{0,\,2}\in C^\infty_{0,\,2}(X,\,\C)$ of $\omega$ represents an $E_2$-cohomology class $\{\rho_\omega^{0,\,2}\}_{E_2}\in E_2^{0,\,2}(X)$. Moreover, $\{\rho_\omega^{0,\,2}\}_{E_2}\in\ker(d_2:E_2^{0,\,2}(X)\to E_2^{2,\,1}(X))$.

   \vspace{1ex}

   (ii)\,  Suppose that $n=3$. Then, the class $\{\rho_\omega^{0,\,2}\}_{E_2}\in E_2^{0,\,2}(X)$ is constant when the Hermitian-symplectic metric $\omega$ varies in a fixed Aeppli cohomology class.

   The class $\{\rho_\omega^{0,\,2}\}_{E_2}\in E_2^{0,\,2}(X)$ will be called the {\bf $E_2$-torsion class} of the Hermitian-symplectic Aeppli class $\{\omega\}_A$.

 \end{Lem-Def}

 \noindent {\it Proof.} (i)\, By construction, the $(0,\,2)$-torsion form $\rho_\omega^{0,\,2}$ has the properties: $$\bar\partial\rho_\omega^{0,\,2} = 0  \hspace{3ex} \mbox{and} \hspace{3ex} \partial\rho_\omega^{0,\,2}\in\mbox{Im}\,\bar\partial \hspace{2ex} (\mbox{since} \hspace{1ex} \partial\rho_\omega^{0,\,2} = -\bar\partial\omega),$$

 \noindent which translate precisely to $\rho_\omega^{0,\,2}$ being $E_2$-closed (see terminology in [Pop19, Proposition 3.1]), namely to $\rho_\omega^{0,\,2}$ representing an $E_2$-cohomology class.

 Moreover, the class $d_2(\{\rho_\omega^{0,\,2}\}_{E_2})\in E_2^{2,\,1}(X)$ is represented by $-\partial\omega$ since $-\omega$ is such that $\bar\partial(-\omega) = \partial\rho_\omega^{0,\,2}$. (See [CFGU97].) However, $\partial\omega$ is $\bar\partial$-exact, so, in particular, $\{\partial\omega\}_{E_2}=0$. We get $$d_2(\{\rho_\omega^{0,\,2}\}_{E_2}) = -\{\partial\omega\}_{E_2}=0\in E_2^{2,\,1}(X).$$

 \vspace{1ex}

 (ii)\, When $n=3$, Corollary \ref{Cor:energy_M-A_mass} tells us that the respective $(2,\,0)$-torsion forms $\rho_\omega^{2,\,0}$ and $\rho_\eta^{2,\,0}:=\rho_{\omega_\eta}^{2,\,0}$ of any two Aeppli-cohomologous Hermitian-symplectic metrics $\omega$ and $\omega_\eta = \omega + \partial\bar\eta + \bar\partial\eta >0$, with $\eta\in C^{\infty}_{1,\,0}(X,\,\C)$, are related by $\rho_\eta^{2,\,0} = \rho_\omega^{2,\,0} + \partial\eta$. Hence, for their conjugates, we get: $$\rho_\eta^{0,\,2} = \rho_\omega^{0,\,2} + \bar\partial\bar\eta, \hspace{3ex} \mbox{so} \hspace{3ex} \{\rho_\eta^{0,\,2}\}_{\bar\partial} = \{\rho_\omega^{0,\,2}\}_{\bar\partial}, \hspace{3ex} \mbox{hence also} \hspace{3ex} \{\rho_\eta^{0,\,2}\}_{E_2} = \{\rho_\omega^{0,\,2}\}_{E_2}\in E_2^{0,\,2}(X).$$  \hfill $\Box$

\vspace{2ex}

Since $\omega$ is K\"ahler if and only if $\rho_\omega^{0,\,2}=0$, we get the following necessary condition for a given Hermitian-symplectic Aeppli class $\{\omega\}_A$ to contain a K\"ahler metric.

\begin{Cor}\label{Cor:necessary-cond_K} Suppose that $n=3$. If a given Hermitian-symplectic Aeppli class $\{\omega\}_A$ contains a  K\"ahler metric, then its $E_2$-torsion class $\{\rho_\omega^{0,\,2}\}_{E_2}\in E_2^{0,\,2}(X)$ vanishes.

  Moreover, the condition $\{\rho_\omega^{0,\,2}\}_{E_2}=0$ in $E_2^{0,\,2}(X)$ is equivalent to $\rho_\omega^{0,\,2}\in\mbox{Im}\,\bar\partial$ for some (hence every) metric $\omega$ lying in $\{\omega\}_A$.

\end{Cor}  

\noindent {\it Proof.} Only the latter statement still needs a proof. The $E_2$-exactness condition on $\rho_\omega^{0,\,2}$ is equivalent to the existence of a $(0,\,1)$-form $\xi$ and of a $(-1,\,2)$-form $\zeta$ such that $\rho_\omega^{0,\,2} = \partial\zeta + \bar\partial\xi$ and $\bar\partial\zeta=0$. (See [Pop19, Proposition 3.1].) However, for bidegree reasons, every $(-1,\,2)$-form $\zeta$ is trivially the zero form, so the $E_2$-exactness condition on $(0,\,2)$-forms is equivalent to the $\bar\partial$-exactness condition.  \hfill $\Box$

\vspace{2ex}

Therefore, we are led to restricting attention to Hermitian-symplectic Aeppli classes of vanishing torsion class on $3$-dimensional manifolds. In this case, $\rho_\omega^{0,\,2}$ is $\bar\partial$-exact and we let \begin{equation}\label{eqn:xi_min-sol_formula}\xi_\omega^{0,\,1} = \Delta^{''-1}\bar\partial^\star\rho_\omega^{0,\,2}\in\mbox{Im}\,\bar\partial^\star\subset C^\infty_{0,\,1}(X,\,\C)\end{equation}

\noindent be the minimal $L^2_\omega$-norm solution of the equation $\bar\partial\xi=\rho_\omega^{0,\,2}$. Our functional $F : {\cal S}_{\{\omega_0\}} \to [0,\,+\infty)$ of Definition \ref{Def:F_energy-functional_H-S} takes the form: \begin{equation}\label{eqn:F_energy-functional_H-S_zero-torsion} F(\omega) = \int\limits_X|\rho_\omega^{2,\,0}|^2_\omega\,dV_\omega = \int\limits_X\rho_\omega^{2,\,0}\wedge\rho_\omega^{0,\,2}\wedge\omega = \int\limits_X\partial\xi_\omega^{1,\,0}\wedge\bar\partial\xi_\omega^{0,\,1}\wedge\omega,\end{equation}

\noindent where $\xi_\omega^{1,\,0}$ is the conjugate of $\xi_\omega^{0,\,1}$.

Meanwhile, using formulae (\ref{eqn:Neumann_formula_del-0}) and (\ref{eqn:xi_min-sol_formula}) and applying the {\it a priori estimate} in various Sobolev norms to the elliptic operator $\Delta''$, we get the following estimates for every $k\in\N$:

\begin{eqnarray}\label{eqn:xi-rho-omega_apriori-estimates}\nonumber||\xi_\omega^{0,\,1}||_{W^{k+2}} & \leq & C_k\,||\bar\partial^\star\rho_\omega^{0,\,2}||_{W^k}\leq A_kC_k\,||\rho_\omega^{0,\,2}||_{W^{k+1}} = A_kC_k\,||\Delta^{'-1}\partial^\star\bar\partial\omega||_{W^{k+1}} \\
  & = & A_kC_k\,||\Delta^{''-1}\bar\partial^\star\partial\omega||_{W^{k+1}} \leq A_kC_k^2\,||\bar\partial^\star\partial\omega||_{W^{k-1}}\leq A_k^2C_k^2\,||\partial\omega||_{W^k},\end{eqnarray}

\noindent where $A_k, C_k>0$ are constants independent of $\omega$. This shows that $\xi_\omega^{0,\,1}$ is small when $\partial\omega$ is small (i.e. $\omega$ is close to being K\"ahler).

\vspace{3ex}

Now, recall that a Hodge theory for the second page $E_2$ of the Fr\"olicher spectral sequence of $X$ was constructed in [Pop16] by means of the introduction of the {\it pseudo-differential Laplace-type operator} $$\widetilde\Delta=\partial p''\partial^\star + \partial^\star p''\partial + \Delta'' : C^\infty_{p,\,q}(X,\,\C)\to C^\infty_{p,\,q}(X,\,\C)$$

\noindent whose kernel in every bidegree $(p,\,q)$ is isomorphic to $E_2^{p,\,q}(X)$, where $p'':C^\infty_{p,\,q}(X,\,\C)\to\ker\Delta''\subset C^\infty_{p,\,q}(X,\,\C)$ is the orthogonal projection onto the $\Delta''$-harmonic space. Thus, every $E_2$-class $\{\alpha\}_{E_2}\in E_2^{p,\,q}$ contains a unique element lying in $$\ker\widetilde\Delta = \ker\bar\partial\cap\ker\bar\partial^\star\cap\ker(p''\partial)\cap\ker(p''\partial^\star).$$

\noindent This construction applies for any fixed Hermitian metric $\omega$ on $X$ with respect to which $\partial^\star$, $\bar\partial^\star$, $p''$, $\Delta''$ and $\widetilde\Delta$ are computed. We will add a subscript $\omega$ when we wish to stress the dependence of the respective operator on $\omega$.

\begin{Lem}\label{Lem:rho_antiharmonicity} Let $X$ be a compact complex manifold such that $\mbox{dim}_\C X=3$. Suppose that there exists a Hermitian-symplectic metric $\omega$ on $X$ such that $\{\rho_\omega^{0,\,2}\}_{E_2}=0$. Then:

  \vspace{2ex}

  (i)\, $\rho_\omega^{0,\,2}\in\ker\bar\partial\cap\ker(p''\partial)\cap\ker(p''\partial^\star)$, hence $\rho_\omega^{0,\,2}\in\ker\widetilde\Delta$ if and only if $\bar\partial^\star\rho_\omega^{0,\,2}=0$;

  \vspace{2ex}

  (ii)\, $\rho_\omega^{0,\,2}\perp\ker\widetilde\Delta$, hence $\rho_\omega^{0,\,2}\in\ker\widetilde\Delta$ if and only if $\rho_\omega^{0,\,2}=0$.

\end{Lem}

\noindent {\it Proof.} (i)\, As noticed in Corollary \ref{Cor:necessary-cond_K}, the hypothesis $\{\rho_\omega^{0,\,2}\}_{E_2}=0$ is equivalent to $\rho_\omega^{0,\,2}\in\mbox{Im}\,\bar\partial$. Hence, $\bar\partial\rho_\omega^{0,\,2}=0$ and $\partial\rho_\omega^{0,\,2}\in\mbox{Im}\,\bar\partial$, so $(p''\partial)\,\rho_\omega^{0,\,2}=0$ as well. Meanwhile, $\rho_\omega^{2,\,0}\in\mbox{Im}\,\bar\partial^\star$ by construction (see Observation \ref{Obs:torsion-formula_dim3}). Taking conjugates, we get $\partial^\star\rho_\omega^{0,\,2}=0$.

\vspace{1ex}

(ii) The following $L^2_\omega$-orthogonal $3$-space decomposition was proved in [Pop16, Lemma 3.3] in every bidegree $(p,\,q)$:

\begin{equation}\label{eqn:3-space-decomp_Delta-tilde}C^{\infty}_{p,\,q}(X,\,\C) = \ker\widetilde{\Delta}\bigoplus\bigg(\mbox{Im}\,\bar\partial + \mbox{Im}\,(\partial_{|\ker\bar\partial})\bigg)\bigoplus\bigg(\mbox{Im}\,(\partial^{\star}\circ p'') +  \mbox{Im}\,\bar\partial^{\star}\bigg).\end{equation}

\noindent In our case, $(p,\,q)=(0,\,2)$, hence $\mbox{Im}\,(\partial_{|\ker\bar\partial}) = \{0\}$ for bidegree reasons. Since $\rho_\omega^{0,\,2}\in\mbox{Im}\,\bar\partial$, we get $\rho_\omega^{0,\,2}\perp\ker\widetilde\Delta$, by orthogonality. \hfill $\Box$

\vspace{3ex}

It follows from the above lemma that the following equivalence holds: $$\rho_\omega^{0,\,2}=0\iff\bar\partial^\star\rho_\omega^{0,\,2}=0.$$

On the other hand, $\bar\partial^\star\rho_\omega^{0,\,2} = \bar\partial^\star\bar\partial\xi_\omega^{0,\,1} = \Delta''\xi_\omega^{0,\,1}$, hence it is tempting to introduce the following

\begin{Def}\label{Def:G-functional} Let $X$ be a compact complex manifold such that $\mbox{dim}_\C X=3$. Suppose that there exists a Hermitian-symplectic metric $\omega_0$ on $X$ such that $\{\rho_{\omega_0}^{0,\,2}\}_{E_2}=0$.

  We define the following {\bf energy functional}: $$G:{\cal S}_{\omega_0}\to[0,\,+\infty),   \hspace{3ex}  G(\omega)=\int\limits_X|\bar\partial^\star_\omega\rho_\omega^{0,\,2}|^2_\omega\,dV_\omega = ||\bar\partial^\star_\omega\rho_\omega^{0,\,2}||^2_\omega.$$

\end{Def}

From the above discussion, we get the identities \begin{equation}\label{eqn:G_identities}G(\omega)=||\bar\partial^\star_\omega\rho_\omega^{0,\,2}||^2_\omega = \langle\langle\Delta''_\omega\rho_\omega^{0,\,2},\,\rho_\omega^{0,\,2}\rangle\rangle_\omega = ||\Delta''_\omega\xi_\omega^{0,\,1}||^2_\omega\end{equation}

\noindent and the equivalences  $$F(\omega)=0 \iff G(\omega)=0 \iff \omega \hspace{1ex}\mbox{is K\"ahler}.$$

\section{Approach via a Monge-Amp\`ere-type equation}\label{section:M-A_eq}

Let $X$ be a compact complex manifold with $\mbox{dim}_\C X=3$. Suppose there exists a Hermitian-symplectic metric $\omega$ on $X$. As usual, let $A = A_{\{\omega\}_A}:= F(\omega) + \mbox{Vol}_\omega(X)>0$ be the generalised volume of the Aeppli cohomology class $\{\omega\}_A\in H^{1,\,1}_A(X,\,\R)$.   

 Let $\gamma$ be an arbitrary Hermitian metric on $X$ and let $c>0$ be the constant defined by \begin{equation}\label{eqn:def_c}\bigg(\int_X\omega\wedge\frac{\gamma^2}{2!}\bigg)^3\bigg\slash\bigg(\int_X\frac{\gamma^3}{3!}\bigg)^2 = 6A/c.\end{equation}

  Now, consider the following Monge-Amp\`ere-type equation: \begin{equation}\label{eqn:M-A_eq}\nonumber (\omega + \partial\bar\eta + \bar\partial\eta)^3 = c\,(\Lambda_\gamma\omega)^3\,\,\frac{\gamma^3}{3!} \hspace{6ex} (\star)\end{equation}

  \noindent whose unknown is $\eta\in C^\infty_{1,\,0}(X,\,\C)$, subject to the positivity condition $\omega + \partial\bar\eta + \bar\partial\eta>0$.

  On the one hand, for any $\eta$ such that $\omega + \partial\bar\eta + \bar\partial\eta>0$, the H\"older inequality with conjugate exponents $p=3$ and $q=3/2$ gives:

  \begin{eqnarray}\label{eqn:Holder_cube}\int\limits_X\sqrt[3]{\frac{(\omega + \partial\bar\eta + \bar\partial\eta)^3}{\gamma^3}}\,\frac{\gamma^3}{3!} \leq \bigg(\int\limits_X\frac{(\omega + \partial\bar\eta + \bar\partial\eta)^3}{3!}\bigg)^{\frac{1}{3}}\,\bigg(\int\limits_X\frac{\gamma^3}{3!}\bigg)^{\frac{2}{3}}.\end{eqnarray}

  On the other hand, if the form $\eta$ solves equation ($\star$), we have \begin{eqnarray}\label{eqn:star_cubic-root_consequence}\int\limits_X\sqrt[3]{\frac{(\omega + \partial\bar\eta + \bar\partial\eta)^3}{\gamma^3}}\,\frac{\gamma^3}{3!} = (c/3!)^{\frac{1}{3}}\, \int\limits_X\Lambda_\gamma\omega\,\frac{\gamma^3}{3!} = (c/3!)^{\frac{1}{3}}\,\int\limits_X\omega\wedge\frac{\gamma^2}{2!}.\end{eqnarray}

  Putting together (\ref{eqn:def_c}), (\ref{eqn:Holder_cube}) and (\ref{eqn:star_cubic-root_consequence}), we get, whenever $\eta$ solves equation ($\star$): \begin{eqnarray}\nonumber A\geq\int\limits_X\frac{(\omega + \partial\bar\eta + \bar\partial\eta)^3}{3!}\geq\frac{1}{(\mbox{Vol}_\gamma(X))^2}\,\bigg(\int\limits_X\sqrt[3]{\frac{(\omega + \partial\bar\eta + \bar\partial\eta)^3}{\gamma^3}}\,\frac{\gamma^3}{3!}\bigg)^3 & = & \frac{c}{3!}\,\frac{\bigg(\int\limits_X\omega\wedge\frac{\gamma^2}{2!}\bigg)^3}{(\mbox{Vol}_\gamma(X))^2 } \\
  \nonumber & = & \frac{c}{3!}\,\frac{6A}{c} = A,\end{eqnarray}

  \noindent where $\mbox{Vol}_\gamma(X) := \int_X\gamma^3/3!$. This implies, thanks to (\ref{eqn:A-invariant}), that $F(\omega_\eta)=0$ in this case, which is equivalent to $\omega_\eta$ being a K\"ahler metric.

  We have thus proved Proposition \ref{Prop:Introd_consequence_MA-eq} stated in the introduction.

\section{Stratification of the Aeppli class}\label{section:stratification}

In $\S.$\ref{section:M-A_eq} we defined the Monge-Amp\`ere type equation $(\star)$ and observed that its solutions yield K\"ahler metrics. Unfortunately, little is known about the solvability of equations of this type. Therefore, in this section  we consider the special case of equation $(\star)$ where the solution $\eta$ is of the shape $\eta = -(i/2)\,\partial\varphi$, with $\varphi:X\to\R$ a $C^\infty$ function. Equation $(\star)$ becomes:  \begin{equation}\label{eqn:M-A_eq_BCc}\nonumber (\omega + i\partial\bar\partial\varphi)^3 = c\,(\Lambda_\gamma\omega)^3\,\,\frac{\gamma^3}{3!}, \hspace{6ex} (\star\star)\end{equation} subject to the condition $\omega + i\partial\bar\partial\varphi>0$, where $\gamma$ is an arbitrary Hermitian metric fixed on $X$, $A$ is the generalised volume of $\{\omega\}_A$ defined in (\ref{eqn:A-invariant}) and the constant $c>0$ is defined in (\ref{eqn:def_c}). The advantage is that we are now dealing with a scalar equation and the existence theory is much more developed in this set-up (c.f. [Che87], [GL09], [TW10]). The drawback is that the perturbation of $\omega$ by $i\partial\bar{\partial}\varphi$ is non-generic within its Aeppli class and this forces us to break $\{\omega\}_A$ into subclasses (to be defined below) and study equation $(\star\star)$ in each subclass.

A conformal rescaling of $\gamma$ by a $C^\infty$ function $f:X\to(0,\,+\infty)$ will change the constant $c>0$ to some constant $c_f>0$ defined by the analogue of (\ref{eqn:def_c}): \begin{equation}\label{eqn:def_c_f}\nonumber\frac{6A}{c_f} = \frac{\bigg(\int\limits_Xf^2\,\omega\wedge\frac{\gamma^2}{2!}\bigg)^3}{\bigg(\int\limits_X\frac{f^3\,\gamma^3}{3!}\bigg)^2} = \frac{\bigg(\int\limits_Xf^2\,(\Lambda_\gamma\omega)\,dV_\gamma\bigg)^3}{\bigg(\int\limits_Xf^3\,dV_\gamma\bigg)^2}\leq\int\limits_X(\Lambda_\gamma\omega)^3\,dV_\gamma,\end{equation} where the last inequality is H\"older's inequality applied to the functions $f^2$ and $\Lambda_\gamma\omega$ with the conjugate exponents $p=3/2$ and $q=3$. This translates to the following eligibility condition for $c_f$: \begin{equation*}c_f\geq\frac{6A}{\int\limits_X(\Lambda_\gamma\omega)^3\,dV_\gamma}.\end{equation*}

\noindent Now, H\"older's inequality is an equality if $f=\Lambda_\gamma\omega$. In particular, if the metric $\gamma>0$ is chosen such that $\Lambda_\gamma\omega\equiv 1$, no conformal rescaling of $\gamma$ is necessary (i.e. we can choose $f\equiv 1$) to get the {\it minimal} constant $$c = \frac{6A}{\int\limits_XdV_\gamma}.$$

\begin{Def}\label{Def:omega-normalised} Let $\omega$ be a fixed Hermitian metric on a compact complex manifold $X$. A Hermitian metric $\gamma$ on $X$ is said to be {\bf $\omega$-normalised} if $\Lambda_\gamma\omega=1$ at every point of $X$.\end{Def}

The following observation is trivial.

\begin{Lem}\label{Lem:conformal-class_omega-normalised} Every conformal class of Hermitian metrics on $X$ {\bf contains} a {\bf unique} $\omega$-normalised representative.

\end{Lem}

\noindent{\it Proof.} Let $\gamma$ be an arbitrary Hermitian metric on $X$. We are looking for $C^\infty$ functions $f:X\to(0,\,+\infty)$ such that $\Lambda_{f\gamma}\omega=1$ on $X$. Since $\Lambda_{f\gamma}\omega=(1/f)\,\Lambda_\gamma\omega$, the only possible choice for $f$ is $f=\Lambda_\gamma\omega$. \hfill $\Box$

\vspace{2ex}

We saw in $\S.$\ref{section:M-A_eq} that if equality is achieved in H\"older's inequality (i.e. if the constant $c>0$ assumes its {\it minimal} value computed above) and if equation $(\star\star)$ is solvable with this minimal constant $c$ on the right, then its solution is a K\"ahler metric. In other words, Proposition \ref{Prop:Introd_consequence_MA-eq} and the above considerations lead to

\begin{Conc}\label{Conc:BC-subclass_after-H-eq} Let $X$ be a compact complex manifold with $\mbox{dim}_\C X=3$. Suppose there exists a {\bf Hermitian-symplectic} metric $\omega$ on $X$ and fix an arbitrary {\bf $\omega$-normalised} Hermitian metric $\gamma$ on $X$. Let $A>0$ be the generalised volume of $\{\omega\}_A$ defined in (\ref{eqn:A-invariant}).

   If there exists a $C^\infty$ solution $\varphi:X\to\R$ of the equation \begin{equation}\label{eqn:M-A_eq_BC}\nonumber \frac{(\omega + i\partial\bar\partial\varphi)^3}{3!} = A\,\frac{dV_\gamma}{\int\limits_XdV_\gamma} \hspace{6ex} (\star\star)\end{equation} such that $\omega_\varphi:=\omega + i\partial\bar\partial\varphi>0$, then $\omega_\varphi$ is a K\"ahler metric lying in the Aeppli cohomology class $\{\omega\}_A$.

\end{Conc}

\subsection{The strata}\label{subsection:strata}

A Hermitian-symplectic metric $\omega$ need not be $d$-closed, but let us still call the affine space $\{\omega\}_{BC}:=\{\omega + i\partial\bar\partial\varphi\,\mid\,\varphi\in C^\infty(X,\,\R)\}$ the {\bf Bott-Chern subclass} (or {\bf stratum}) of $\omega$. It is a subspace of the Aeppli class $\{\omega\}_A$ of $\omega$. Similarly, by analogy with the open convex subset $${\cal S}_{\{\omega\}}:=\{\omega'>0\,\mid\,\omega'\in\{\omega\}_A\}\subset\{\omega\}_A\cap C^\infty_{1,\,1}(X,\,\R)$$ of metrics in the Aeppli class of a given H-S metric $\omega$, we define the open convex subset $${\cal D}_{[\omega']}:=\{\omega''>0\,\mid\,\omega''\in\{\omega'\}_{BC}\}\subset\{\omega'\}_{BC}$$ of metrics in the Bott-Chern subclass of a given H-S metric $\omega'$.

If we fix a Hermitian-symplectic metric $\omega$, we can {\bf partition} ${\cal S}_{\{\omega\}}$ as \begin{equation}\label{eqn:stratification}{\cal S}_{\{\omega\}}=\bigcup\limits_{j\in J}{\cal D}_{[\omega_j]},\end{equation} where $(\omega_j)_{j\in J}$ is a system of representatives of the Bott-Chern subclasses ${\cal D}_{[\omega']}$ when $\omega'$ ranges over ${\cal S}_{\{\omega\}}$. Moreover, for every $j\in J$, let $\gamma_j$ be an $\omega_j$-normalised Hermitian metric on $X$ and let us consider the equation: \begin{equation*}\frac{(\omega_j + i\partial\bar\partial\varphi)^3}{3!} = A\,\frac{dV_{\gamma_j}}{\int\limits_XdV_{\gamma_j}} \hspace{6ex} (\star\star_j)\end{equation*} such that $\omega_j + i\partial\bar\partial\varphi>0$. (No other condition is imposed at this point on $\gamma_j$.) 

By the Tosatti-Weinkove theorem [TW10, Corollary 1], there exists a unique constant $b_j>0$ such that the equation \begin{equation*}\frac{(\omega_j + i\partial\bar\partial\varphi)^3}{3!} = b_jA\,\frac{dV_{\gamma_j}}{\int\limits_XdV_{\gamma_j}}\hspace{6ex} (\star\star\star_j),\end{equation*} subject to the extra condition $\omega_j + i\partial\bar\partial\varphi>0$, is solvable. Integrating and using the inequality $\int_X(\omega_j + i\partial\bar\partial\varphi)^3/3!\leq A$, which follows from (\ref{eqn:A-invariant}), we get: $$b_j\leq 1,  \hspace{3ex} j\in J.$$

From this and from Conclusion \ref{Conc:BC-subclass_after-H-eq}, we infer Proposition \ref{Prop:Introd_b_j1} stated in the Introduction.

\vspace{2ex}

The next observation is that, within Bott-Chern subclasses of Hermitian-symplectic metrics that contain a {\it Gauduchon} metric, the volume remains {\it constant} and all the metrics are {\it Gauduchon}. These Bott-Chern subclasses will be called {\it Gauduchon strata}. 

\begin{Lem}\label{Lem:BC-subclass_volume-Gauduchon} Let $\mbox{dim}_\C X=3$. Suppose that a metric $\omega$ on $X$ is both {\bf SKT} and {\bf Gauduchon}. Then, for every $\varphi:X\longrightarrow\R$, we have $$(a)\,\int\limits_X(\omega + i\partial\bar\partial\varphi)^3 = \int\limits_X\omega^3  \hspace{3ex} \mbox{and} \hspace{3ex} (b)\,\,\,\partial\bar\partial(\omega + i\partial\bar\partial\varphi)^2 = 0.$$    \end{Lem}

\noindent {\it Proof.} (a)\, Straightforward calculations give: \begin{eqnarray*}\int\limits_X(\omega + i\partial\bar\partial\varphi)^3 & = & \int\limits_X\omega^3 + 3\int\limits_X\omega^2\wedge i\partial\bar\partial\varphi + 3\int\limits_X\omega\wedge(i\partial\bar\partial\varphi)^2 + \int\limits_X(i\partial\bar\partial\varphi)^3  \\
  & = & \int\limits_X\omega^3 + 3i\int\limits_X\varphi\,\partial\bar\partial\omega^2 - 3\int\limits_X\varphi\,\partial\bar\partial\omega\wedge\partial\bar\partial\varphi + \int\limits_X\partial(i\bar\partial\varphi\wedge(i\partial\bar\partial\varphi)^2) = \int\limits_X\omega^3,\end{eqnarray*}

\noindent where $\partial\bar\partial\omega=0$ since $\omega$ is SKT, while $\partial\bar\partial\omega^2=0$ since $\omega$ is Gauduchon.

\vspace{1ex}

(b)\, Straightforward calculations give: \begin{eqnarray*}\partial\bar\partial(\omega + i\partial\bar\partial\varphi)^2 = \partial\bar\partial\omega^2 + 2\,\partial\bar\partial\omega\wedge(i\partial\bar\partial\varphi) + \partial\bar\partial(i\partial\bar\partial\varphi)^2 =0,\end{eqnarray*} since $\partial\bar\partial\omega=0$ and $\partial\bar\partial\omega^2=0$ for the same reasons as in (a).  \hfill $\Box$

\vspace{2ex}

\subsection{Volume comparison within a Bott-Chern stratum}\label{subsection:volume_BC-stratum}

We have seen that for any SKT metric $\omega$ on a $3$-dimensional compact complex manifold $X$, we have: \begin{eqnarray}\label{eqn:vol_BC_comparison1}\nonumber\int\limits_X\frac{(\omega + i\partial\bar\partial\varphi)^3}{3!} & = & \int\limits_X\frac{\omega^3}{3!} + \int\limits_X\omega^2\wedge\frac{i}{2}\partial\bar\partial\varphi = \int\limits_X\frac{\omega^3}{3!} + \int\limits_X \varphi\,\frac{i}{2}\partial\bar\partial\omega^2 \\
  & = & \int\limits_X\frac{\omega^3}{3!} + \int\limits_X \varphi\,i\partial\omega\wedge\bar\partial\omega,\end{eqnarray} where the SKT hypothesis on $\omega$ is used only to get the last identity. Thus, to understand the variation of $\mbox{Vol}_{\omega_\varphi}(X)=\int_X(\omega + i\partial\bar\partial\varphi)^3/3!$ when $\varphi$ ranges over the $C^\infty$ real-valued functions on $X$ such that $\omega + i\partial\bar\partial\varphi>0$, the following observation will come in handy.

\begin{Lem}\label{Lem:almost_lck-bal_formula} Let $X$ be a $3$-dimensional complex manifold and let $\omega$ be an arbitrary Hermitian metric on $X$. If  $$\partial\omega = (\partial\omega)_{prim} + \alpha^{1,\,0}\wedge\omega$$ is the Lefschetz decomposition of $\partial\omega$ into a primitive part (w.r.t. $\omega$) and a part divisible by $\omega$, with $\alpha^{1,\,0}\in C^\infty_{1,\,0}(X,\,\C)$, then \begin{eqnarray}\label{eqn:almost_lck-bal_formula}i\partial\omega\wedge\bar\partial\omega = \bigg(|\alpha^{1,\,0}\wedge\omega|^2_\omega - |(\partial\omega)_{prim}|^2_\omega\bigg)\,dV_\omega.\end{eqnarray}

\end{Lem}

\noindent {\it Proof.} From the Lefschetz decomposition, we get: \begin{eqnarray*}i\partial\omega\wedge\bar\partial\omega = i(\partial\omega)_{prim}\wedge(\bar\partial\omega)_{prim} + i\alpha^{1,\,0}\wedge\alpha^{0,\,1}\wedge\omega^2,\end{eqnarray*} where $\alpha^{0,\,1}=\overline{\alpha^{1,\,0}}$. This is because $(\partial\omega)_{prim}\wedge\omega=0$ and $(\bar\partial\omega)_{prim}\wedge\omega=0$. Indeed, $(\partial\omega)_{prim}$ and $(\bar\partial\omega)_{prim}$ are primitive $3$-forms on a $3$-dimensional complex manifold, so they lie in the kernel of $\omega\wedge\cdot$.

From the general formula (\ref{eqn:prim-form-star-formula-gen}), we get: $$(\bar\partial\omega)_{prim} = i\,\star(\bar\partial\omega)_{prim} \hspace{3ex} \mbox{and} \hspace{3ex} \alpha^{0,\,1}\wedge\frac{\omega^2}{2!} = -i\,\star\alpha^{0,\,1},$$ where $\star$ is the Hodge star operator induced by $\omega$. Hence, \begin{eqnarray*}i(\partial\omega)_{prim}\wedge(\bar\partial\omega)_{prim} & = & -(\partial\omega)_{prim}\wedge\star(\bar\partial\omega)_{prim} = - |(\partial\omega)_{prim}|^2\,dV_\omega \\
i\alpha^{1,\,0}\wedge\alpha^{0,\,1}\wedge\omega^2 & = & 2\,\alpha^{1,\,0}\wedge\star\alpha^{0,\,1} = 2\,|\alpha^{1,\,0}|^2_\omega\,dV_\omega.\end{eqnarray*}

Formula (\ref{eqn:almost_lck-bal_formula}) follows from these computations after further noticing that \begin{eqnarray*}|\alpha^{1,\,0}\wedge\omega|^2_\omega & = & \langle\alpha^{1,\,0}\wedge\omega,\,\alpha^{1,\,0}\wedge\omega\rangle_\omega = \langle\Lambda_\omega(\alpha^{1,\,0}\wedge\omega),\,\alpha^{1,\,0}\rangle_\omega \\
  & = & \langle[\Lambda_\omega,\,\omega\wedge\cdot](\alpha^{1,\,0}),\,\alpha^{1,\,0}\rangle_\omega = 2\,|\alpha^{1,\,0}|^2_\omega,\end{eqnarray*} where the last identity follows from the well-known formula $[\Lambda_\omega,\,\omega\wedge\cdot] = (n-k)\,\mbox{Id}$ on $k$-forms on an $n$-dimensional complex manifold. (In our case, $k=1$ and $n=3$.)  \hfill $\Box$

\vspace{3ex}

Notice that, in the setting of Lemma \ref{Lem:almost_lck-bal_formula}, $\omega$ is {\it balanced} if and only if $\partial\omega = (\partial\omega)_{prim}$, while $\omega$ is {\it lck} (i.e. {\it locally conformally K\"ahler}) if and only if $\partial\omega = \alpha^{1,\,0}\wedge\omega$. This accounts for the terminology used in the next

\begin{Cor}\label{Cor:almost_lck-bal} Let $X$ be a $3$-dimensional compact complex manifold equipped with an {\bf SKT metric} $\omega$.

  \vspace{1ex}

  (a)\, If $|\alpha^{1,\,0}\wedge\omega|_\omega\geq|(\partial\omega)_{prim}|_\omega$ at every point of $X$ (we will say in this case that $\omega$ is {\bf almost lck}), then $\omega$ is {\bf Gauduchon} and $|\alpha^{1,\,0}\wedge\omega|_\omega=|(\partial\omega)_{prim}|_\omega$ at every point of $X$.

  \vspace{1ex}

  (b)\, If $|(\partial\omega)_{prim}|_\omega\geq|\alpha^{1,\,0}\wedge\omega|_\omega$ at every point of $X$ (we will say in this case that $\omega$ is {\bf almost balanced}), then $\omega$ is {\bf Gauduchon} and $|\alpha^{1,\,0}\wedge\omega|_\omega=|(\partial\omega)_{prim}|_\omega$ at every point of $X$.

  \vspace{1ex}

  (c)\, $\omega$ is {\bf almost lck} $\iff$ $\omega$ is {\bf almost balanced} $\iff$ $\omega$ is {\bf Gauduchon} $\iff$ $|\alpha^{1,\,0}\wedge\omega|_\omega=|(\partial\omega)_{prim}|_\omega$ at every point of $X$.

\end{Cor}  

\noindent {\it Proof.} The SKT assumption on $\omega$ implies that $i\partial\omega\wedge\bar\partial\omega = \frac{i}{2}\,\partial\bar\partial\omega^2$. Integrating this identity and using the Stokes theorem and formula (\ref{eqn:almost_lck-bal_formula}), we get: \begin{eqnarray}\label{eqn:integrals-equal_Lefschetz}\int\limits_X|\alpha^{1,\,0}\wedge\omega|^2_\omega\,dV_\omega = \int\limits_X|(\partial\omega)_{prim}|_\omega^2\,dV_\omega.\end{eqnarray}

Therefore, if $|(\partial\omega)_{prim}|^2_\omega-|\alpha^{1,\,0}\wedge\omega|^2_\omega$ has constant sign on $X$, it must vanish identically. This is equivalent to $\frac{i}{2}\,\partial\bar\partial\omega^2$ vanishing identically, hence to $\omega$ being Gauduchon.  \hfill $\Box$

\vspace{3ex}

Based on these observations, let us introduce the following

\begin{Not}\label{Not:puddles} For any SKT metric $\omega$ on a $3$-dimensional compact complex manifold $X$, we put:  \begin{eqnarray*}U_\omega: &=&\bigg\{x\in X\,\bigm| \,|\alpha^{1,\,0}\wedge\omega|_\omega(x) < |(\partial\omega)_{prim}|_\omega(x)\bigg\}, \\
    V_\omega: &=&\bigg\{x\in X\,\bigm| \,|\alpha^{1,\,0}\wedge\omega|_\omega(x) > |(\partial\omega)_{prim}|_\omega(x)\bigg\}, \\
  Z_\omega: &=&\bigg\{x\in X\,\bigm| \,|\alpha^{1,\,0}\wedge\omega|_\omega(x) = |(\partial\omega)_{prim}|_\omega(x)\bigg\}.\end{eqnarray*}

\end{Not}  

Clearly, $U_\omega$ and $V_\omega$ are open subsets of $X$, while $Z_\omega$ is closed. The three of them form a partition of $X$. Moreover, Corollary \ref{Cor:almost_lck-bal} ensures that $\omega$ is Gauduchon if and only
if $U_\omega = V_\omega = \emptyset$. This happens if and only if either $U_\omega = \emptyset$ or $V_\omega = \emptyset$.

\vspace{2ex}

Returning to the variation of the volume of $\omega_\varphi:=\omega + i\partial\bar\partial\varphi$, we now observe a stark contrast between the non-Gauduchon strata dealt with below and the Gauduchon ones treated in Lemma \ref{Lem:BC-subclass_volume-Gauduchon}.

\begin{Lem}\label{Lem:BC-subclass_volume-nonGauduchon} Let $X$ be a $3$-dimensional compact complex manifold. Suppose that $\omega$ is an {\bf SKT non-Gauduchon} metric on $X$.
  Then, the map $$\bigg\{\varphi\in C^\infty(X)\,\big|\omega + i\partial\bar\partial\varphi>0\bigg\}\ni\varphi\longmapsto\int\limits_X\frac{(\omega + i\partial\bar\partial\varphi)^3}{3!}:={Vol}_{\omega_\varphi}(X)\in(0,\,+\infty)$$ does not achieve any local extremum.

\end{Lem}

\noindent {\it Proof.} Suppose this map achieves, say, a local maximum at some metric $\omega_0=\omega + i\partial\bar\partial\varphi_0>0$. Without loss of generality, we may assume that
$\omega_0=\omega$ (and $\varphi_0\equiv 1$). Since $\omega$ is not Gauduchon, both $U_\omega$ and $V_\omega$ are not empty. Thanks to (\ref{eqn:vol_BC_comparison1}) and (\ref{eqn:almost_lck-bal_formula}), the
local maximality of $\omega$ translates to \begin{eqnarray}\label{eqn:local-max_ineq}\int\limits_X\varphi\,|\alpha^{1,\,0}\wedge\omega|^2_\omega\,dV_\omega \leq \int\limits_X\varphi\,|(\partial\omega)_{prim}|_\omega^2\,dV_\omega\end{eqnarray}
for every $\varphi\in C^\infty(X,\,\R)$ such that $\omega + i\partial\bar\partial\varphi>0$ and $\varphi$ is close enough to $\varphi_0\equiv 1$ in $C^2$ norm.

  Now, (\ref{eqn:integrals-equal_Lefschetz}) translates to \begin{eqnarray}\label{eqn:local-max_ineq_1}\nonumber\int\limits_{U_\omega}(|\alpha^{1,\,0}\wedge\omega|^2_\omega - |(\partial\omega)_{prim}|_\omega^2)\,dV_\omega & + & \int\limits_{V_\omega}(|\alpha^{1,\,0}\wedge\omega|^2_\omega - |(\partial\omega)_{prim}|_\omega^2)\,dV_\omega \\
    & + & \int\limits_{Z_\omega}(|\alpha^{1,\,0}\wedge\omega|^2_\omega - |(\partial\omega)_{prim}|_\omega^2)\,dV_\omega = 0.\end{eqnarray}

  Thus, if we can find a $\varphi\in C^\infty(X,\,\R)$ sufficiently close to $\varphi_0\equiv 1$ in $C^2$ norm (this will also imply that $\omega + i\partial\bar\partial\varphi>0$) such that $$\varphi\equiv 1 \hspace{1ex}\mbox{on}\hspace{1ex} U_\omega\cup Z_\omega, \hspace{3ex} \varphi\equiv 1+\varepsilon \hspace{1ex}\mbox{on}\hspace{1ex} V'_\omega\Subset V_\omega, \hspace{3ex} \mbox{and}\hspace{3ex} 1\leq\varphi\leq 1+\varepsilon \hspace{1ex}\mbox{on}\hspace{1ex} V_\omega\setminus V'_\omega,$$ for some constant $\varepsilon>0$, where $V'_\omega$ is a pregiven relatively compact open subset of $V_\omega$, we will have $$\int\limits_{V_\omega}(\varphi - 1)\,(|\alpha^{1,\,0}\wedge\omega|^2_\omega-|(\partial\omega)_{prim}|_\omega^2)\,dV_\omega>0.$$ Thanks to (\ref{eqn:local-max_ineq_1}), this will imply that \begin{eqnarray*}\int\limits_X\varphi\,|\alpha^{1,\,0}\wedge\omega|^2_\omega\,dV_\omega > \int\limits_X\varphi\,|(\partial\omega)_{prim}|_\omega^2\,dV_\omega,\end{eqnarray*} which will contradict (\ref{eqn:local-max_ineq}).

  Now, if $\varepsilon>0$ is chosen small enough, it is obvious that a function $\varphi\in C^\infty(X,\,\R)$ with the above properties exists. \hfill $\Box$

  \vspace{3ex}

  Summing up, the volume of $\omega_\varphi:=\omega + i\partial\bar\partial\varphi$ is {\it constant} on the {\it Gauduchon strata} (if any), while it {\it achieves no local extremum} on the {\it non-Gauduchon strata}.  

 \section{Cohomological interpretations of the generalised volume}\label{Gen-volume_cohom}

Before turning to these interpretations of our invariant $A$ in $\S.$\ref{subsection:cohom_A} and $\S.$\ref{subsection:minimal-completion}, we first display it in the context of Hermitian-symplectic and strongly Gauduchon metrics in $\S.$\ref{subsection:sG_H-S}.

\subsection{sG metrics induced by H-S metrics}\label{subsection:sG_H-S}

From Proposition \ref{Prop:H-S_sG}, we infer the following construction. Let $\mbox{dim}_\C X=3$. With any Hermitian-symplectic metric $\omega$ on $X$, we uniquely associate the $C^\infty$ positive definite $(2,\,2)$-form \begin{equation}\label{eqn:Omega-omega_def}\Omega_\omega:=\omega^2 + 2\rho_\omega^{2,\,0}\wedge\rho_\omega^{0,\,2},\end{equation} where $\rho_\omega^{2,\,0}$ is the $(2,\,0)$-{\it torsion form} of $\omega$ and $\rho_\omega^{0,\,2}=\overline{\rho_\omega^{2,\,0}}$. As is well known (see e.g. [Mic83]), there exists a unique positive definite $(1,\,1)$-form $\gamma_\omega$ such that $$\gamma_\omega^2 =\Omega_\omega.$$ By construction and the proof of Proposition \ref{Prop:H-S_sG}, $\gamma_\omega$ is a {\it strongly Gauduchon} metric on $X$ that will be called the {\bf sG metric associated with $\omega$}. Of course, $\gamma_\omega = \omega$ if and only if $\omega$ is K\"ahler. Since $\gamma_\omega^2$ and $\Omega_\omega$ determine each other uniquely, we will often identify them. In particular, we will also refer to $\Omega_\omega$ as the {\it sG metric associated with $\omega$}. We get: $$\frac{1}{3!}\,\Omega_\omega\wedge\omega = \frac{1}{3!}\,\omega^3 + \frac{1}{3}\,|\rho_\omega^{2,\,0}|_\omega^2\,dV_\omega.$$ Hence, \begin{equation}\label{eqn:Omega_omega-omega}\frac{1}{6}\,\int\limits_X\Omega_\omega\wedge\omega = \frac{2}{3}\,\mbox{Vol}_\omega(X) + \frac{1}{3}\,A,\end{equation} where $A=\mbox{Vol}_\omega(X) + F(\omega)>0$ is the generalised volume of the H-S Aeppli class $\{\omega\}_A$.

Thus, the problem of maximising $\mbox{Vol}_\omega(X)$ when $\omega$ ranges over the metrics in  $\{\omega\}_A$ is equivalent to maximising the quantity $\int_X\Omega_\omega\wedge\omega$.

\subsection{The first cohomological interpretation of the generalised volume $A$}\label{subsection:cohom_A}

We first observe that the Aeppli cohomology class of $\Omega_\omega$ depends only on the Aeppli class of $\omega$.

\begin{Lem}\label{Lem:Aeppli-class_Omega_omega} Suppose that $\mbox{dim}_\C X=3$. For any {\bf Aeppli cohomologous} Hermitian-symplectic metrics $\omega$ and $\omega_\eta = \omega + \partial\bar\eta + \bar\partial\eta$ on $X$, with $\eta\in C^\infty_{1,\,0}(X,\,\C)$, the associated sG metrics $\Omega_\omega$ and $\Omega_{\omega_\eta}$ are again {\bf Aeppli cohomologous}.

  Specifically, we have: 
  \begin{eqnarray}\label{eqn:Aeppli-class_Omega_omega}\nonumber\Omega_{\omega_\eta} - \Omega_\omega & = & \partial(\bar\eta\wedge\partial\bar\eta) + \bar\partial(\eta\wedge\bar\partial\eta) + 2\,\partial(\bar\eta\wedge\bar\partial\eta) + 2\,\bar\partial(\partial\eta\wedge\bar\eta) \\
    & + & 2\,\partial(\eta\wedge\rho^{0,\,2}_\omega) + 2\,\bar\partial(\bar\eta\wedge\rho^{2,\,0}_\omega) + 2\,\partial(\bar\eta\wedge\omega) + 2\,\bar\partial(\eta\wedge\omega),\end{eqnarray} so $\Omega_{\omega_\eta} - \Omega_\omega\in\mbox{Im}\,\partial + \mbox{Im}\,\bar\partial$.

\end{Lem}

\noindent {\it Proof.} We know from Corollary \ref{Cor:energy_M-A_mass} that the $(2,\,0)$-torsion forms of $\omega_\eta$ and $\omega$ are related by $\rho^{2,\,0}_\eta = \rho^{2,\,0}_\omega + \partial\eta$. We get: \begin{eqnarray}\nonumber\Omega_{\omega_\eta} & = & \omega_\eta^2 + 2\,\rho^{2,\,0}_{\omega_\eta}\wedge\rho^{0,\,2}_{\omega_\eta} = (\omega + \partial\bar\eta + \bar\partial\eta)^2 + 2\,(\rho^{2,\,0}_\omega + \partial\eta)\wedge(\rho^{0,\,2}_\omega + \bar\partial\bar\eta) \\
  \nonumber & = & \omega^2 + (\partial\bar\eta + \bar\partial\eta)^2 + 2\,\omega\wedge(\partial\bar\eta + \bar\partial\eta) + 2\,\rho^{2,\,0}_\omega\wedge\rho^{0,\,2}_\omega + 2\,\partial\eta\wedge\bar\partial\bar\eta + 2\,\rho^{2,\,0}_\omega\wedge\bar\partial\bar\eta + 2\,\partial\eta\wedge\rho^{0,\,2}_\omega \\
  \nonumber & = & \Omega_\omega + \partial(\bar\eta\wedge\partial\bar\eta) + \bar\partial(\eta\wedge\bar\partial\eta) + 2\,\partial(\bar\eta\wedge\bar\partial\eta) + 2\,\bar\eta\wedge\partial\bar\partial\eta \\
  \nonumber & + & 2\,\partial\bar\eta\wedge\omega + 2\,\bar\partial\eta\wedge\omega + 2\,\bar\partial(\partial\eta\wedge\bar\eta) + 2\,\partial\bar\partial\eta\wedge\bar\eta \\
   \nonumber & + & 2\,\bar\partial(\bar\eta\wedge\rho^{2,\,0}_\omega) + 2\,\partial(\bar\eta\wedge\omega) - 2\,\partial\bar\eta\wedge\omega + 2\,\partial(\eta\wedge\rho^{0,\,2}_\omega) + 2\,\bar\partial(\eta\wedge\omega) - 2\,\bar\partial\eta\wedge\omega.\end{eqnarray} This proves (\ref{eqn:Aeppli-class_Omega_omega}) since all the terms that are neither in $\mbox{Im}\,\partial$ nor in $\mbox{Im}\,\bar\partial$ reoccur with the opposite sign and cancel. \hfill $\Box$

\vspace{3ex}

We will now use some notions introduced in [PSU20]. We will need the following

\begin{Def}(Definition 3.1 in [PSU20b])\label{Def:E_rE_r-bar} Let $X$ be a compact complex manifold. Fix any $r\geq 1$.

\vspace{1ex}

(i)\, A form $\alpha\in C^\infty_{p,\,q}(X)$ is said to be {\bf $E_r\overline{E}_r$-closed} if $\partial\bar\partial\alpha=0$ and if there exist smooth forms $\eta_1,\dots , \eta_{r-1}$ and $\rho_1,\dots , \rho_{r-1}$ such that the following two towers of $r-1$ equations are satisfied: \begin{align*}\label{eqn:towers_E_rE_r-bar-closedness} \partial\alpha & = \bar\partial\eta_1 &    \bar\partial\alpha & = \partial\rho_1 & \\
    \partial\eta_1 & = \bar\partial\eta_2 &   \bar\partial\rho_1 & = \partial\rho_2 & \\
     \vdots & & \\
     \partial\eta_{r-2} & =  \bar\partial\eta_{r-1},  &  \bar\partial\rho_{r-2} & = \partial\rho_{r-1}. &\end{align*}

\vspace{1ex}

(ii)\, A form $\alpha\in C^\infty_{p,\,q}(X)$ is said to be {\bf $E_r\overline{E}_r$-exact} if there exist smooth forms $\zeta, \xi, \eta$ such that \begin{equation}\label{eqn:main-eq_E_rE_r-bar-exactness}\alpha = \partial\zeta + \partial\bar\partial\xi + \bar\partial\eta\end{equation}

\noindent and such that $\zeta$ and $\eta$ further satisfy the following conditions. There exist smooth forms $v_{r-3},\dots , v_0$ and $u_{r-3},\dots , u_0$ such that the following two towers of $r-1$ equations are satisfied: \begin{align*}\label{eqn:towers_E_rE_r-bar-exactness} \bar\partial\zeta & = \partial v_{r-3} &    \partial\eta & = \bar\partial u_{r-3} & \\
     \bar\partial v_{r-3} & = \partial v_{r-4} &   \partial u_{r-3} & = \bar\partial u_{r-4} & \\
     \vdots & & \\
     \bar\partial v_0 & =  0,  &  \partial u_0 & = 0. &\end{align*}  

\end{Def}

\vspace{2ex}

Note that when $r=1$, the two towers in (i) are empty, so the $E_1\overline{E}_1$-closedness condition is equivalent to $\partial\bar\partial$-closedness. Meanwhile, when $r\geq 2$, the two towers in (i) imply that $\partial\bar\partial\alpha = 0$. As for (ii), when $r=1$, the $E_1\overline{E}_1$-exactness condition is equivalent to $\partial\bar\partial$-exactness.

We will also need the following

\begin{Def}(Definition 3.4 in [PSU20b])\label{Def:E_r-BC_E_r-A} Let $X$ be a compact complex manifold. Fix $r\in\N^\star$ and a bidegree $(p,\,q)$.

\vspace{1ex}

(i)\, The {\bf $E_r$-Bott-Chern} cohomology group of bidegree $(p,\,q)$ of $X$ is defined as the following quotient complex vector space: \[E_{r,\,BC}^{p,\,q}(X):=\frac{\{\alpha\in C^\infty_{p,\,q}(X)\,\mid\,d\alpha=0\}}{\{\alpha\in C^\infty_{p,\,q}(X)\,\mid\,\alpha\hspace{1ex}\mbox{is}\hspace{1ex} E_r\overline{E}_r\mbox{-exact}\}}.\]

\vspace{1ex}

(ii)\, The {\bf $E_r$-Aeppli} cohomology group of bidegree $(p,\,q)$ of $X$ is defined as the following quotient complex vector space: \[E_{r,\,A}^{p,\,q}(X):=\frac{\{\alpha\in C^\infty_{p,\,q}(X)\,\mid\,\alpha\hspace{1ex}\mbox{is}\hspace{1ex} E_r\overline{E}_r-\mbox{closed}\}}{\{\alpha\in C^\infty_{p,\,q}(X)\,\mid\,\alpha\in\mbox{Im}\,\partial + \mbox{Im}\,\bar\partial\}}.\]

\end{Def}

We will also need the following

\begin{Lem}(Proposition 6.2 in [PSU20b])\label{Lem:sG-HS_E_2A} Let $X$ be a compact complex manifold with $\mbox{dim}_\C X=n$ and let $\omega$ be a Hermitian metric on $X$.

  \vspace{1ex}

  (i)\, The metric $\omega$ is {\bf strongly Gauduchon (sG)} if and only if $\omega^{n-1}$ is {\bf $E_2\overline{E}_2$-closed}.

  \vspace{1ex}

  (ii)\, The metric $\omega$ is {\bf Hermitian-symplectic (H-S)} if and only if $\omega$ is {\bf $E_3\overline{E}_3$-closed}.

  If $n=3$, $\omega$ is {\bf Hermitian-symplectic (H-S)} if and only if $\omega$ is {\bf $E_2\overline{E}_2$-closed}.

\end{Lem}




\vspace{3ex}

In our case, the consequence of Proposition \ref{Prop:H-S_sG} and of Lemmas \ref{Lem:Aeppli-class_Omega_omega} and \ref{Lem:sG-HS_E_2A} is the following

\begin{Lem}\label{Lem:Omega_omega_same-E2A-class} Let $X$ be a compact complex manifold with $\mbox{dim}_\C X=3$. For any {\bf Aeppli cohomologous} Hermitian-symplectic metrics $\omega$ and $\omega_\eta = \omega + \partial\bar\eta + \bar\partial\eta$ on $X$, with $\eta\in C^\infty_{1,\,0}(X,\,\C)$, the corresponding sG metrics $\Omega_\omega$ and $\Omega_{\omega_\eta}$ represent the {\bf same $E_2$-Aeppli class}: $$\{\Omega_{\omega_\eta}\}_{E_{2,\,A}} = \{\Omega_\omega\}_{E_{2,\,A}}\in E_{2,\,A}^{2,\,2}(X).$$

\end{Lem}

\noindent {\it Proof.} We know from Proposition \ref{Prop:H-S_sG} that $\Omega_\omega$ and $\Omega_{\omega_\eta}$ are sG metrics, so by (i) of Lemma \ref{Lem:sG-HS_E_2A} they represent $E_2$-Aeppli classes. Meanwhile,  by Lemma \ref{Lem:Aeppli-class_Omega_omega}, $\Omega_\omega$ and $\Omega_{\omega_\eta}$ are Aeppli cohomologous, hence also $E_2$-Aeppli cohomologous. \hfill $\Box$

\vspace{3ex}

We will now use the main notion introduced in [PSU20a].

\begin{Def}(Theorem and Definition 1.2 in [PSU20a])\label{Def:page-r-ddbar_def} Fix an arbitrary $r\in\N^\star$. A $n$-dimensional compact complex manifold $X$ is said to be a {\bf page-$(r-1)$-$\partial\bar\partial$-manifold} if $X$ has the {\bf $E_r$-Hodge Decomposition} property in the following sense.

 For every bidegree $(p,\,q)$, every class $\{\alpha^{p,\,q}\}_{E_r}\in E_r^{p,\,q}(X)$ can be represented by a {\bf $d$-closed $(p,\,q)$-form} and for every $k$, the linear map
 
$$\bigoplus_{p+q=k}E_r^{p,\,q}(X)\ni\sum\limits_{p+q=k}\{\alpha^{p,\,q}\}_{E_r}\mapsto\bigg\{\sum\limits_{p+q=k}\alpha^{p,\,q}\bigg\}_{DR}\in H^k_{DR}(X,\,\C)$$

\noindent is {\bf well-defined} by means of $d$-closed pure-type representatives and {\bf bijective}.

\end{Def}

\vspace{3ex}

We will also need the following result from [PSU20].

\begin{The}(Theorem 3.53 in [PSU20])\label{The:page_r-1_ddbar_prop-B_BC-A} Let $X$ be a compact complex manifold with $\mbox{dim}_\C X=n$. Fix an arbitrary integer $r\geq 1$. The following statements are {\bf equivalent}.

  \vspace{1ex}

  (a)\, $X$ is a {\bf page-$(r-1)$-$\partial\bar\partial$-manifold}.

  \vspace{1ex}

  (b)\, For all $p,q\in\{0,\dots , n\}$, the canonical linear maps $T_r^{p,\,q}:E_{r,\,BC}^{p,\,q}(X)\longrightarrow E_r^{p,\,q}(X)$ and $S_r^{p,\,q}:E_r^{p,\,q}(X)\longrightarrow E_{r,\,A}^{p,\,q}(X)$ induced by the identity are {\bf isomorphisms}.

\end{The}

\vspace{3ex}

In our case, as a consequence of Theorem \ref{The:page_r-1_ddbar_prop-B_BC-A}, we get a {\it unique lift} $\mathfrak{c}_\omega\in E_{2,\,BC}^{2,\,2}(X)$ of $\{\Omega_\omega\}_{E_{2,\,A}}\in E_{2,\,A}^{2,\,2}(X)$ under the appropriate assumption on $X$. 

\begin{Cor}\label{Cor:E_2BC_lifts} Let $X$ be a {\bf page-$1$-$\partial\bar\partial$-manifold} with $\mbox{dim}_\C X=3$. For any {\bf Aeppli cohomologous} Hermitian-symplectic metrics $\omega$ and $\omega_\eta = \omega + \partial\bar\eta + \bar\partial\eta$ on $X$, with $\eta\in C^\infty_{1,\,0}(X,\,\C)$, there exists a unique {\bf $E_2$-Bott-Chern class} $\mathfrak{c}_\omega\in E_{2,\,BC}^{2,\,2}(X)$ such that $$(S_2^{2,\,2}\circ T_2^{2,\,2})(\mathfrak{c}_\omega) =  \{\Omega_{\omega_\eta}\}_{E_{2,\,A}} = \{\Omega_\omega\}_{E_{2,\,A}}\in E_{2,\,A}^{2,\,2}(X),$$ where $\Omega_\omega$ and $\Omega_{\omega_\eta}$ are the sG metrics associated with $\omega$, resp. $\omega_\eta$.

  In particular, the {\bf $E_2$-Bott-Chern class} $\mathfrak{c}_\omega\in E_{2,\,BC}^{2,\,2}(X)$ depends only on the {\bf $E_2$-Aeppli class} $\{\omega\}_{E_2,\,A}\in E_{2,\,A}^{1,\,1}(X)$.

\end{Cor}

\vspace{2ex}

We can now state and prove the main result of this subsection. It will use the {\bf duality} between the $E_r$-Bott-Chern cohomology of any bidegree $(p,\,q)$ and the $E_r$-Aeppli cohomology of the complementary bidegree $(n-p,\,n-q)$ proved as Theorem 3.11 in [PSU20b]. In our case, $n=3$, $r=2$ and $(p,\,q)=(2,\,2)$.

\begin{The}\label{The:cohom_A} Let $X$ be a {\bf page-$1$-$\partial\bar\partial$-manifold} with $\mbox{dim}_\C X=3$. Suppose there exists a \newline Hermitian-symplectic metric $\omega$ on $X$ whose {\bf $E_2$-torsion class vanishes} (i.e. $\{\rho^{0,\,2}_\omega\}_{E_2} = 0\in E_2^{0,\,2}(X)$).

  Then, the generalised volume $A=F(\omega) + \mbox{Vol}_\omega(X)$ of $\{\omega\}_A$ is given as the following intersection number in cohomology: \begin{equation}\label{eqn:cohom_A}A=\frac{1}{6}\,\mathfrak{c}_\omega.\{\omega\}_{E_2,\,A}.\end{equation}

\end{The}

\noindent {\it Proof.} $\bullet$ We will first construct a smooth $d$-closed $(2,\,2)$-form $\widetilde\Omega_\omega$ that represents the $E_2$-Bott-Chern class $\mathfrak{c}_\omega\in E_{2,\,BC}^{2,\,2}(X)$ in the most economical way possible. We will proceed in two stages that correspond to lifting the $E_2$-Aeppli class $\{\Omega_\omega\}_{E_2,\,A}\in E_{2,\,A}^{2,\,2}(X)$ to $E_2^{2,\,2}(X)$ under the isomorphism $S_2^{2,\,2}:E_2^{2,\,2}(X)\longrightarrow E_{2,\,A}^{2,\,2}(X)$ induced by the identity, respectively to lifting the resulting $E_2$-class in $E_2^{2,\,2}(X)$ to $E_{2,\,BC}^{2,\,2}(X)$ under the isomorphism $T_2^{2,\,2}:E_{2,\,BC}^{2,\,2}(X)\longrightarrow E_2^{2,\,2}(X)$ induced by the identity.

\vspace{1ex}

{\it Stage $1$.} To lift $\{\Omega_\omega\}_{E_2,\,A}\in E_{2,\,A}^{2,\,2}(X)$ to $E_2^{2,\,2}(X)$ under the isomorphism $S_2^{2,\,2}:E_2^{2,\,2}(X)\longrightarrow E_{2,\,A}^{2,\,2}(X)$, we need to find a $(2,\,2)$-form $\Gamma^{2,\,2}$ such that $\Gamma^{2,\,2}\in\mbox{Im}\,\partial + \mbox{Im}\,\bar\partial$ (because we need $\Omega_\omega + \Gamma^{2,\,2}$ to represent the same $E_{2,\,A}$-class as the original $\Omega_\omega$) and such that $$\bar\partial(\Omega_\omega + \Gamma^{2,\,2}) = 0 \hspace{3ex} \mbox{and} \hspace{3ex} \partial(\Omega_\omega + \Gamma^{2,\,2})\in\mbox{Im}\,\bar\partial,$$ (because we need $\Omega_\omega + \Gamma^{2,\,2}$ to represent an $E_2$-class). The last two conditions are equivalent to \begin{eqnarray}\label{eqn:cohom_A_proof_1}\bar\partial\Gamma^{2,\,2} = -\bar\partial\Omega_\omega \hspace{3ex} \mbox{and} \hspace{3ex} \partial\Gamma^{2,\,2}\in\mbox{Im}\,\bar\partial,\end{eqnarray} because, for the last condition, we already have $\partial\Omega_\omega\in\mbox{Im}\,\bar\partial$ by the sG property of $\Omega_\omega$.

If we denote by ${\cal Z}_{2\bar{2}}^{2,\,2}$ the space of smooth $E_2\overline{E}_2$-closed $(2,\,2)$-forms on $X$, we have $\Omega_\omega\in{\cal Z}_{2\bar{2}}^{2,\,2}$, hence $-\bar\partial\Omega_\omega\in\bar\partial({\cal Z}_{2\bar{2}}^{2,\,2})$. On the other hand, part (ii) of Lemma 3.52 in [PSU20] ensures that $\bar\partial({\cal Z}_{2\bar{2}}^{2,\,2})\subset\mbox{Im}\,(\partial\bar\partial)$ because $X$ is a {\it page-$1$-$\partial\bar\partial$-manifold}. (Actually, this inclusion is equivalent to the surjectivity of the map $S_2^{2,\,2}$.) We conclude that $-\bar\partial\Omega_\omega\in\mbox{Im}\,(\partial\bar\partial)$, so the equation \begin{equation}\label{eqn:u12_def}\partial\bar\partial u^{1,\,2} = \bar\partial\Omega_\omega\end{equation} admits solutions $u^{1,\,2}\in C^\infty_{1,\,2}(X,\,\C)$. Let $u^{1,\,2}_\omega$ be the minimal $L^2_\omega$-norm such solution and put $\Gamma^{2,\,2}_\omega:=\partial u^{1,\,2}_\omega$.

Thus, $\Gamma^{2,\,2}_\omega=\partial u^{1,\,2}_\omega$ satisfies conditions (\ref{eqn:cohom_A_proof_1}) and $\Gamma^{2,\,2}_\omega\in\mbox{Im}\,\partial\subset\mbox{Im}\,\partial + \mbox{Im}\,\bar\partial$. So, we have got the minimal lift $\{\Omega_\omega + \partial u^{1,\,2}_\omega\}_{E_2}$ of $\{\Omega_\omega\}_{E_2,\,A}\in E_{2,\,A}^{2,\,2}(X)$ to $E_2^{2,\,2}(X)$, i.e. $$S_2^{2,\,2}(\{\Omega_\omega + \partial u^{1,\,2}_\omega\}_{E_2}) = \{\Omega_\omega\}_{E_2,\,A}.$$

\vspace{1ex}

{\it Stage $2$.} To lift $\{\Omega_\omega + \partial u^{1,\,2}_\omega\}_{E_2}\in E_2^{2,\,2}(X)$ to $E_{2,\,BC}^{2,\,2}(X)$ under the isomorphism \newline $T_2^{2,\,2}:E_{2,\,BC}^{2,\,2}(X)\longrightarrow E_2^{2,\,2}(X)$, we need to find a $(2,\,2)$-form $V^{2,\,2}$ such that $V^{2,\,2}\in\partial(\ker\bar\partial) + \mbox{Im}\,\bar\partial$ (because we need $\Omega_\omega + \partial u^{1,\,2}_\omega + V^{2,\,2}$ to represent the same $E_2$-class as $\Omega_\omega + \partial u^{1,\,2}_\omega$) and such that $$\partial(\Omega_\omega + \partial u^{1,\,2}_\omega + V^{2,\,2}) = 0 \hspace{3ex} \mbox{and} \hspace{3ex} \bar\partial(\Omega_\omega + \partial u^{1,\,2}_\omega + V^{2,\,2}) = 0,$$ (because we need $\Omega_\omega + \partial u^{1,\,2}_\omega + V^{2,\,2}$ to represent an $E_{2,\,BC}$-class). The last two conditions are equivalent to \begin{eqnarray}\label{eqn:cohom_A_proof_2}\partial V^{2,\,2} = -\partial(\Omega_\omega + \partial u^{1,\,2}_\omega) \hspace{3ex} \mbox{and} \hspace{3ex} \bar\partial V^{2,\,2} = 0,\end{eqnarray} because, for the last condition, we already have $\bar\partial(\Omega_\omega + \partial u^{1,\,2}_\omega) = 0$.

Now, $\Omega_\omega + \partial u^{1,\,2}_\omega\in{\cal Z}_{2\bar{2}}^{2,\,2}$, hence $-\partial(\Omega_\omega + \partial u^{1,\,2}_\omega)\in\partial({\cal Z}_{2\bar{2}}^{2,\,2})\subset\mbox{Im}\,(\partial\bar\partial)$, the last inclusion being a consequence of part (ii) of Lemma 3.52 in [PSU20] and of $X$ being a {\it page-$1$-$\partial\bar\partial$-manifold}.  (Actually, this inclusion is equivalent to the surjectivity of the map $T_2^{2,\,2}$.)

We conclude that $-\partial(\Omega_\omega + \partial u^{1,\,2}_\omega)\in\mbox{Im}\,(\partial\bar\partial)$, so the equation \begin{equation}\label{eqn:u21_def}\partial\bar\partial u^{2,\,1} = -\partial(\Omega_\omega + \partial u^{1,\,2}_\omega)\end{equation} admits solutions $u^{2,\,1}\in C^\infty_{2,\,1}(X,\,\C)$. Let $u^{2,\,1}_\omega$ be the minimal $L^2_\omega$-norm such solution and put $V^{2,\,2}_\omega:=\bar\partial u^{2,\,1}_\omega$. Clearly, $u^{2,\,1}_\omega = \overline{u^{1,\,2}_\omega}$ since the equations whose minimal solutions are $u^{2,\,1}_\omega$ and $u^{1,\,2}_\omega$ are conjugated to each other.

Thus, $V^{2,\,2}_\omega=\bar\partial u^{2,\,1}_\omega$ satisfies conditions (\ref{eqn:cohom_A_proof_2}) and $V^{2,\,2}_\omega\in\partial(\ker\bar\partial) + \mbox{Im}\,\bar\partial$.

\vspace{2ex}

The upshot of the above construction is that the form $$\widetilde\Omega_\omega:=\Omega_\omega + \partial u^{1,\,2}_\omega + \bar\partial u^{2,\,1}_\omega\in C^\infty_{2,\,2}(X,\,\C)$$ is the minimal completion of $\Omega_\omega$ to a {\it $d$-closed} pure-type form of bidegree $(2,\,2)$. Moreover, the class $\{\widetilde\Omega_\omega\}_{E_{2,\,BC}}\in E_{2,\,BC}^{2,\,2}(X)$ has the property that $$(S_2^{2,\,2}\circ T_2^{2,\,2})(\{\widetilde\Omega_\omega\}_{E_{2,\, BC}}) = \{\Omega_\omega\}_{E_2,\,A}.$$ Hence, $\{\widetilde\Omega_\omega\}_{E_{2,\,BC}} = \mathfrak{c}_\omega$ since the map $S_2^{2,\,2}\circ T_2^{2,\,2}:E_{2,\, BC}^{2,\,2}(X)\longrightarrow E_{2,\, A}^{2,\,2}(X)$ is bijective (thanks to the {\it page-$1$-$\partial\bar\partial$-assumption} on $X$).  

\vspace{2ex}

$\bullet$ We will now use the representative $\widetilde\Omega_\omega$ of the class $\mathfrak{c}_\omega$ to relate the intersection number in (\ref{eqn:cohom_A}) to the generalised volume of $\{\omega\}_A$. We have: \begin{eqnarray*}\mathfrak{c}_\omega.\{\omega\}_{E_2,\,A} = \int\limits_X\widetilde\Omega_\omega\wedge\omega = \int\limits_X\Omega_\omega\wedge\omega + \int\limits_X\partial u^{1,\,2}_\omega\wedge\omega + \int\limits_X\bar\partial u^{2,\,1}_\omega\wedge\omega.\end{eqnarray*}

Since $\rho^{2,\,0}_\omega = \partial\xi^{1,\,0}_\omega$ thanks to the hypothesis $\{\rho^{0,\,2}_\omega\}_{E_2} = 0\in E_2^{0,\,2}(X)$ (see Corollary \ref{Cor:necessary-cond_K} and the minimal choice (\ref{eqn:xi_min-sol_formula}) of $\xi^{0,\,1}_\omega=\overline{\xi^{1,\,0}_\omega}$), we get: \begin{eqnarray*}\int\limits_X\partial u^{1,\,2}_\omega\wedge\omega & = & \int\limits_Xu^{1,\,2}_\omega\wedge\partial\omega = -\int\limits_Xu^{1,\,2}_\omega\wedge\bar\partial\rho^{2,\,0}_\omega = \int\limits_Xu^{1,\,2}_\omega\wedge\partial\bar\partial\xi^{1,\,0}_\omega \\
  & = & \int\limits_X\partial\bar\partial u^{1,\,2}_\omega\wedge\xi^{1,\,0}_\omega \stackrel{(a)}{=} \int\limits_X\bar\partial\Omega_\omega\wedge\xi^{1,\,0}_\omega \stackrel{(b)}{=} -2\,\int\limits_X\partial(\rho^{0,\,2}_\omega\wedge\omega)\wedge\xi^{1,\,0}_\omega \\
  & = & 2\,\int\limits_X\rho^{0,\,2}_\omega\wedge\omega\wedge\partial\xi^{1,\,0}_\omega = 2\,\int\limits_X\rho^{2,\,0}_\omega\wedge\rho^{0,\,2}_\omega\wedge\omega = 2\,||\rho^{2,\,0}_\omega||^2_\omega = 2F(\omega),\end{eqnarray*} where (a) follows from (\ref{eqn:u12_def}) and (b) follows from the formula \begin{equation}\label{eqn:dbar_Omega-omega_formula}\bar\partial\Omega_\omega = -2\,\partial(\rho^{0,\,2}_\omega\wedge\omega)\end{equation} which in turn follows at once from (\ref{eqn:Omega-omega_def}).

By conjugation, we also have $\int_X\bar\partial u^{2,\,1}_\omega\wedge\omega = 2F(\omega)$. Putting the various pieces of information together, we get: \begin{eqnarray*}\mathfrak{c}_\omega.\{\omega\}_{E_2,\,A} =  \int\limits_X\Omega_\omega\wedge\omega + 4F(\omega) = 4\,\mbox{Vol}_\omega(X) + 2A + 4F(\omega) = 6A,\end{eqnarray*} where the second identity follows from (\ref{eqn:Omega_omega-omega}).

The proof of Theorem \ref{The:cohom_A} is complete.  \hfill $\Box$



\subsection{The second cohomological interpretation of the generalised volume}\label{subsection:minimal-completion}

We will now work in the general case (i.e. without the extra assumptions made in Theorem \ref{The:cohom_A}). The result will show, yet again, that the generalised volume $A = A_{\{\omega\}_A}>0$ of a Hermitian-symplectic Aeppli class $\{\omega\}_A$ is a natural analogue in this more general context of the volume of a K\"ahler class. 

\begin{Def}\label{Def:minimal-completion} Let $X$ be a compact complex manifold with $\mbox{dim}_\C X=3$. For any Hermitian-symplectic metric $\omega$ on $X$, the $d$-closed real $2$-form $$\widetilde\omega = \rho_\omega^{2,\,0} + \omega + \rho_\omega^{0,\,2}$$ is called the {\bf minimal completion} of $\omega$, where $\rho_\omega^{2,\,0}$, resp. $\rho_\omega^{0,\,2}$, is the $(2,\,0)$-torsion form, resp.  the $(0,\,2)$-torsion form, of $\omega$.\end{Def}

We will now notice the following consequence of Corollary \ref{Cor:energy_M-A_mass}  It gives a new cohomological interpretation of the generalised volume $A = A_{\{\omega\}_A}>0$.

\begin{Prop}\label{Prop:same-DR-class} Let $X$ be a compact complex Hermitian-symplectic manifold of dimension $n=3$.

  \vspace{1ex}

  (a)\, For any Hermitian-symplectic metric $\omega$ on $X$, its minimal completion $2$-form $\widetilde\omega$ has \newline the property: \begin{equation}\label{eqn:min-comp_integral}\int\limits_X\frac{\widetilde\omega^3}{3!} = \mbox{Vol}_\omega(X) + F(\omega) = A_{\{\omega\}_A}.\end{equation}

  \vspace{1ex}

 (b)\, For any Aeppli-cohomologous Hermitian-symplectic metrics $\omega$ and $\omega_\eta$ \begin{equation}\omega_\eta = \omega + \partial\bar\eta + \bar\partial\eta >0  \hspace{3ex} (\mbox{where} \hspace{1ex} \eta\in C^{\infty}_{1,\,0}(X,\,\C)),\end{equation}

  \noindent the respective minimal completion $2$-forms $\widetilde\omega_\eta$ and $\widetilde\omega$ lie in the same De Rham cohomology class.

  \vspace{1ex}

  Thus, $A_{\{\omega\}_A} = \{\widetilde\omega\}_{DR}^3/3!$.

 \end{Prop}

\noindent {\it Proof.} (a)\, Using (\ref{eqn:Omega_omega-omega}) for identity (a) below and the above notation, we get: \begin{eqnarray*}\int\limits_X \widetilde\omega^3 & = & \int\limits_X\widetilde\omega^2\wedge(\rho_\omega^{2,\,0} + \omega + \rho_\omega^{0,\,2}) = \int\limits_X\Omega_\omega\wedge(\rho_\omega^{2,\,0} + \omega + \rho_\omega^{0,\,2})\\
&+& 2\,\int\limits_X\rho_\omega^{2,\,0}\wedge\omega\wedge\rho_\omega^{0,\,2} + 2\,\int\limits_X\rho_\omega^{0,\,2}\wedge\omega\wedge\rho_\omega^{2,\,0}  \\
  & = & \int\limits_X\Omega_\omega\wedge\omega + 4 F(\omega) \stackrel{(a)}{=} 4\mbox{Vol}_\omega(X) + 2\mbox{Vol}_\omega(X)+ 2F(\omega) + 4 F(\omega) = 6A.\end{eqnarray*}

\vspace{1ex}

(b)\, We know from Corollary \ref{Cor:energy_M-A_mass} that the $(2,\,0)$-torsion forms of $\omega_\eta$ and $\omega$ are related by $\rho^{2,\,0}_\eta = \rho^{2,\,0}_\omega + \partial\eta$. We get: \begin{eqnarray*}\widetilde\omega_\eta = \rho^{2,\,0}_\eta + \omega_\eta + \rho^{0,\,2}_\eta = \widetilde\omega + d(\eta + \bar\eta).\end{eqnarray*}

This proves the contention.  \hfill $\Box$

\vspace{3ex}

\noindent {\bf References.} \\


\vspace{1ex}

\noindent [Che87]\, P. Cherrier --- {\it \'Equations de Monge-Amp\`ere sur les vari\'et\'es hermitiennes compactes} --- Bull. Sc. Math. (2) {\bf 111} (1987), 343-385.

\vspace{1ex}

\noindent [Dem97]\, J.-P. Demailly --- {\it Complex Analytic and Algebraic Geometry} --- http://www-fourier.ujf-grenoble.fr/~demailly/books.html

\vspace{1ex}

\noindent [Don06]\, S. K. Donaldson --- {\it Two-forms on Four-manifolds and Elliptic Equations} --- Inspired by S. S. Chern, 153–172, Nankai Tracts Math.,11, World Sci.Publ., Hackensack, NJ, 2006.

\vspace{1ex}

\noindent [EFV12]\, N. Enrietti, A. Fino, L. Vezzoni --- {\it Tamed Symplectic Forms and Strong K\"ahler with Torsion Metrics} --- J. Symplectic Geom. {\bf 10}, No. 2 (2012) 203-223. 

\vspace{1ex}

\noindent [Gau77a]\, P. Gauduchon --- {\it Le th\'eor\`eme de l'excentricit\'e nulle} --- C.R. Acad. Sc. Paris, S\'erie A, t. {\bf 285} (1977), 387-390.

\vspace{1ex}

\noindent [Gau77b]\, P. Gauduchon --- {\it Fibr\'es hermitiens \`a endomorphisme de Ricci non n\'egatif} --- Bull. Soc. Math. France {\bf 105} (1977) 113-140.

\vspace{1ex}

\noindent [GL09]\, B. Guan, Q. Li --- {\it Complex Monge-Amp\`ere Equations on Hermitian Manifolds} --- arXiv e-print DG 0906.3548v1.

\vspace{1ex}

\noindent [IP13]\, S. Ivanov, G. Papadopoulos --- {\it Vanishing Theorems on $(l/k)$-strong K\"ahler Manifolds with Torsion} ---  Adv.Math. {\bf 237} (2013) 147-164.

\vspace{1ex}

\noindent[HL83]\, R. Harvey, H. B. Lawson --- {\it An intrinsic characterization of K\"ahler manifolds} --- Invent. Math. {\bf 74}, (1983), 169-198.

\vspace{1ex}

\noindent [KS60]\, K. Kodaira, D.C. Spencer --- {\it On Deformations of Complex Analytic Structures, III. Stability Theorems for Complex Structures} --- Ann. of Math. {\bf 71}, no.1 (1960), 43-76.

\vspace{1ex}

\noindent [Lam99]\, A. Lamari --- {\it Courants k\"ahl\'eriens et surfaces compactes} --- Ann. Inst. Fourier {\bf 49}, no. 1 (1999), 263-285.

\vspace{1ex}

\noindent [LZ09]\, T.-J. Li, W. Zhang --- {\it Comparing Tamed and Compatible Symplectic Cones and Cohomological Properties of Almost Complex Manifolds} --- Comm. Anal. Geom. {\bf 17}, no. 4 (2009), 651–683.

\vspace{1ex}

\noindent [Mic83]\, M. L. Michelsohn --- {\it On the Existence of Special Metrics in Complex Geometry} --- Acta Math. {\bf 143} (1983) 261-295.

\vspace{1ex}

\noindent [Pop13]\, D. Popovici --- {\it Deformation Limits of Projective Manifolds: Hodge Numbers and Strongly Gauduchon Metrics} --- Invent. Math. {\bf 194} (2013), 515-534.

\vspace{1ex}

\noindent [Pop15]\, D. Popovici --- {\it Aeppli Cohomology Classes Associated with Gauduchon Metrics on Compact Complex Manifolds} --- Bull. Soc. Math. France {\bf 143}, no. 3 (2015), 1-37.

\vspace{1ex}

\noindent [Pop16]\, D. Popovici --- {\it Degeneration at $E_2$ of Certain Spectral Sequences} ---  Internat. J. of Math. {\bf 27}, no. 13 (2016), DOI: 10.1142/S0129167X16501111.

\vspace{1ex}

\noindent [PSU20]\, D. Popovici, J. Stelzig, L. Ugarte --- {\it Some Aspects of Higher-Page Non-K\"ahler Hodge Theory} --- arXiv e-print AG 2001.02313v1.

\vspace{1ex}

\noindent [PSU20a]\, D. Popovici, J. Stelzig, L. Ugarte --- {\it Higher-Page Hodge Theory of Compact Complex Manifolds} --- arXiv e-print AG 2001.02313v2.

\vspace{1ex}

\noindent [PSU20b]\, D. Popovici, J. Stelzig, L. Ugarte --- {\it Higher-Page Bott-Chern and Aeppli Cohomologies and Applications} ---  arXiv e-print AG 2007.03320v1.

\vspace{1ex}

\noindent [Sch07]\, M. Schweitzer --- {\it Autour de la cohomologie de Bott-Chern} --- arXiv e-print math.AG/0709.3528v1.

\vspace{1ex}

\noindent [ST10]\, J. Streets, G. Tian --- {\it A Parabolic Flow of Pluriclosed Metrics} --- Int. Math. Res. Notices, {\bf 16} (2010), 3101-3133.

\vspace{1ex}

\noindent[Sul76]\, D. Sullivan ---{\it Cycles for the dynamical study of foliated manifolds and complex manifolds}--- Invent. Math. {\bf 36} (1976), 225-255. 

\vspace{1ex}

\noindent [TW10]\, V. Tosatti, B. Weinkove --- {\it The Complex Monge-Amp\`ere Equation on Compact Hermitian Manifolds} --- J. Amer. Math. Soc. 23 (2010), no. 4, 1187-1195.

\vspace{1ex}

\noindent [Ver14]\, M. Verbitsky --- {\it Rational Curves and Special Metrics on Twistor Spaces} --- Geometry and Topology {\bf 18} (2014), 897–909.

\vspace{1ex}

\noindent [Voi02]\, C. Voisin --- {\it Hodge Theory and Complex Algebraic Geometry. I.} --- Cambridge Studies in Advanced Mathematics, 76, Cambridge University Press, Cambridge, 2002.

\vspace{1ex}

\noindent [YZZ19]\, S.-T. Yau, Q. Zhao, F. Zheng --- {\it On Strominger K\"ahler-like Manifolds with Degenerate Torsion} --- arXiv e-print DG 1908.05322v2

\vspace{6ex}

\noindent Department of Mathematics and Computer Science     \hfill Institut de Math\'ematiques de Toulouse,

\noindent Jagiellonian University     \hfill  Universit\'e Paul Sabatier,

\noindent 30-409 Krak\'ow, Ul. Lojasiewicza 6, Poland    \hfill  118 route de Narbonne, 31062 Toulouse, France

\noindent  Email: Slawomir.Dinew@im.uj.edu.pl                \hfill     Email: popovici@math.univ-toulouse.fr

\end{document}